\documentclass{article}
\usepackage{graphicx}
\usepackage{geometry}
\usepackage{amsmath, amsthm, amssymb}
\usepackage{amsfonts}

\usepackage{titlesec}
\usepackage{enumitem}
\usepackage{mathrsfs}

\usepackage{xcolor}
\usepackage{tcolorbox}

\usepackage{marvosym}
\usepackage{soul}
\usepackage{dsfont}
\usepackage{tikz-cd}
\usepackage{natbib}

\definecolor{deepblue}{RGB}{0,51,153}
\usepackage[
    colorlinks=true,
    linkcolor=deepblue,
    citecolor=deepblue,
    urlcolor=deepblue
]{hyperref}

\newtheorem{proposition}{Proposition}
\newtheorem{theorem}{Theorem}

\newtheorem{lemma}{Lemma}

\newtheorem{assumption}{Assumption}

\theoremstyle{plain}
\newtheorem{definition}{Definition}
\newtheorem{remark}{Remark}[section]

\newcommand{\tk}{\mathrm{k}}
\newcommand{\bx}{\boldsymbol{\mathrm{x}}}

\newcommand{\br}{{\boldsymbol{\mathrm{r}}}}

\newcommand{\ud}{\,\mathrm{d}}

\newcommand{\JKO}{\mathcal{J}\mathrm{KO}}

\newenvironment{coloredbox}[2]
{\begin{center}\begin{tcolorbox}[enhanced,breakable,colframe=#1,title=#2,fonttitle=\bfseries,width=0.95\textwidth]}
{\end{tcolorbox}\end{center}}

\definecolor{cblindCB}{HTML}{1F449C}
\definecolor{vertForet}{rgb}{0.03,0.5,0.05}

\titleformat{\section}
{\normalfont\Large\bfseries} 
{\thesection}                
{1em}                        
{}                           

\title{Convergence of Implicit Particle Methods for Wasserstein Gradient Flows}

\author{Yunhong Deng\thanks{School of Mathematics, University of Minnesota,
Minneapolis, MN 55455. (\href{mailto:deng0335@umn.edu}{deng0335@umn.edu}, \href{mailto:liwang@umn.edu}{liwang@umn.edu})}
\and Wilfrid Gangbo\thanks{Department of Mathematics, University of California Los Angeles, Los Angeles, CA 90095.\\ (\href{mailto:wgangbo@math.ucla.edu}{wgangbo@math.ucla.edu})}
\and Li Wang\footnotemark[1]}

\date{\today}

\begin{document}
\maketitle

\begin{abstract}
In this paper, we establish the convergence of implicit particle methods for Wasserstein gradient flows as the number of particles tends to infinity, with the time step held fixed. The schemes combine the Jordan--Kinderlehrer--Otto (JKO) variational time discretization with particle approximations of the evolving probability measure. The main challenge lies in approximating the internal energy, which depends nonlinearly on the probability density and cannot, in general, be evaluated directly on empirical measures. We consider three approaches: a one-dimensional density reconstruction based on the spacings between ordered particles; a blob method using kernel-based density reconstruction; and a transport-map formulation in which the density is recovered from the pushforward relation through the change-of-variables formula. The latter two approaches apply in arbitrary spatial dimensions without requiring the spatial domain to be compact. For the blob method, we also establish simultaneous convergence as the particle number tends to infinity and the regularization parameter vanishes, under the scaling $N\varepsilon_N^d\to\infty$. Our analysis combines $\Gamma$-convergence and compactness arguments to establish convergence of the discrete variational problems and their minimizers. These results complement classical dynamical mean-field theory by providing a variational particle-to-continuum limit for implicit JKO schemes involving nonlinear internal energies.
\end{abstract}



\section{Introduction}
Many phenomena in the biological and physical sciences, including chemotactic migration of bacterial colonies \cite{keller1970initiation, carrillo2019population}, collective swarming of animal groups \cite{carrillo2010particle}, and granular flows of inelastic particles \cite{benedetto1997kinetic, toscani2000one}, can be modeled by collective motion of agents or particles along the negative gradient of a density-dependent potential field $\phi(\mu)$:
\begin{equation}\label{0}
    \frac{\partial \mu}{\partial t}
    - \operatorname{div}\big(\mu \nabla \phi(\mu)\big) = 0,
    \quad \text{on } \mathbb{R}^{d} \times (0, \infty),
\end{equation}
where $\mu(x,t)$ denotes the population or mass density. In many models, the potential takes the form
\begin{equation*}
    \phi(\mu) := U^{\prime}(\mu) + V + W \ast \mu.
\end{equation*}
Here, $U : \mathbb{R}_{+} \to \mathbb{R}$ is the internal energy density, accounting for local effects such as particle collisions, stochastic noise, or volume exclusion; $V : \mathbb{R}^{d} \to \mathbb{R}$ is a confinement potential encoding environmental effects; and $W : \mathbb{R}^{d} \to \mathbb{R}$ is a symmetric interaction potential describing nonlocal attraction or repulsion among agents or particles. Specific assumptions and choices of $U$, $V$, and $W$ are given in Assumption~\ref{assumption} below.

A central feature of~\eqref{0} is its interpretation as a gradient flow in the Wasserstein space $(P_{2}(\mathbb{R}^{d}), W_{2})$ \cite{JKO, otto2001geometry}, where $P_{2}(\mathbb{R}^{d})$ denotes the space of probability measures with finite second moments and $W_{2}$ is the $2$-Wasserstein distance. The corresponding free energy is
\begin{equation}\label{2}
    \mathcal{F}(\mu)
    = \mathcal{F}_{U}(\mu)
    + \mathcal{F}_{V}(\mu)
    + \mathcal{F}_{W}(\mu),
\end{equation}
where the confinement and interaction energies are
\begin{equation} \label{FVW}
    \mathcal{F}_{V}(\mu)
    := \int_{\mathbb{R}^{d}} V(x) \ud \mu(x),
    \qquad
    \mathcal{F}_{W}(\mu)
    := \frac{1}{2}
    \int_{\mathbb{R}^{d} \times \mathbb{R}^{d}}
    W(x-y) \ud \mu(x) \ud \mu(y).
\end{equation}
For an absolutely continuous measure $\mu$, identified with its density, the internal energy is given by
\begin{equation} \label{FU}
    \mathcal{F}_{U}(\mu)
    := \int_{\mathbb{R}^{d}} U(\mu(x)) \ud x.
\end{equation}
This variational perspective leads to the minimizing movements, known in this setting as the Jordan--Kinderlehrer--Otto (JKO) scheme \cite{JKO}, which provides a time-discrete approximation of~\eqref{0}.
\begin{definition}[JKO scheme]\label{def:JKO}
Given a time step $\tau > 0$ and an initial measure
$\mu_{0} \in P_{2}(\mathbb{R}^{d})$, the JKO scheme constructs
a sequence $(\mu_{\tk})_{\tk \geq 0}$ through the minimization problems
\begin{equation}\label{1}
    \begin{aligned}
        \mu_{\tk+1}
        &\in
        \operatorname*{argmin}_{\mu \in P_{2}(\mathbb{R}^{d})}
        \JKO(\mu,\mu_{\tk}):=   \mathcal{F}(\mu)
        + \frac{1}{2\tau} W_{2}(\mu,\mu_{\tk})^{2}.
    \end{aligned}
\end{equation}
The associated piecewise constant interpolation
$\mu_{\tau} : [0,\infty) \to P_{2}(\mathbb{R}^{d})$ is defined by
$\mu_{\tau}(0)=\mu_{0}$ and
\begin{equation*}
    \mu_{\tau}(t) := \mu_{\tk},
    \qquad
    (\tk-1)\tau < t \leq \tk\tau,
    \quad \tk \geq 1.
\end{equation*}
Here, the $2$-Wasserstein distance is given by
\begin{equation*}
    W_{2}(\mu,\nu)
    :=
    \left[
        \min_{\pi \in \Pi(\mu,\nu)}
        \int_{\mathbb{R}^{d} \times \mathbb{R}^{d}}
        |x-y|^{2} \ud \pi(x,y)
    \right]^{1/2},
\end{equation*}
where $\Pi(\mu,\nu)$ denotes the set of probability measures on
$\mathbb{R}^{d} \times \mathbb{R}^{d}$ with marginals $\mu$ and $\nu$.
\end{definition}

The Wasserstein gradient flow framework has led to a rich analytical theory of~\eqref{0} (see for example the overview \cite{carrillo2019aggregation}). Under suitable assumptions, it provides a variational approach to constructing global-in-time, measure-valued weak solutions to~\eqref{0} \cite{blanchet2008convergence, carrillo2011global}. Central to this approach is the convergence of the time-discrete JKO approximations $\mu_{\tau}(t)$ to solutions $\mu(t)$ of~\eqref{0} as $\tau\to0$. It has been established for a range of models, including the Fokker--Planck equation \cite{JKO} and the Keller--Segel equation \cite{blanchet2008convergence}; see also \cite{ambrosio2008gradient} for a general theory of gradient flows in metric spaces. Beyond constructing solutions, the framework connects quantitative long-time behavior to the geometry of the free energy $\mathcal{F}$ through functional inequalities, including log-Sobolev-type inequalities \cite{carrillo2003kinetic, carrillo2004wasserstein, carrillo2006contractions}.

The same variational structure motivated a broad class of structure-preserving numerical methods for~\eqref{0} \cite{carrillo2021lagrangian, carrillo2022primal}. It supplies a natural framework of constructing high-order time-stepping schemes \cite{westdickenberg2010variational, matthes2019variational, chu2025guaranteeing, han2023high} and attracted increasing interest for its potential to alleviate the computational challenges associated with high-dimensional problems \cite{lee2024deep, feng2026learn}.


In this work, we study the convergence of spatially discrete versions of the JKO scheme~\eqref{1} as the spatial resolution is refined. We focus on particle discretizations where the solution of the JKO scheme is approximated using the average of Dirac delta distributions (or empirical measures). The principal technical difficulty lies in constructing a consistent discrete approximation of the internal energy $\mathcal{F}_{U}$, since this energy is generally not well defined on empirical measures. We consider three representative approaches. The first is a one-dimensional particle discretization based on the spacing between neighboring particles \cite{blanchet2008convergence, westdickenberg2010variational, di2015rigorous}. The second is the blob method, in which the empirical measure is regularized before evaluating the internal energy \cite{carrillo2019blob, craig2023blob}. The third is based on transport-map representations, where the density and hence the internal energy are expressed through the deformation of the initial measure \cite{JungeMatthesOsberger2017, lee2024deep}.

\subsection{Particle methods}
We consider a particle approximation of the JKO scheme~\eqref{1}, in which the measure $\mu_{\tk}$ is approximated by an empirical measure. More precisely, let $\{x_{\tk,\ell}\}_{\ell=1}^{N} \subset \mathbb{R}^{d}$ denote $N$ particle locations representing $\mu_{\tk}$, and define
\begin{equation} \label{mu_kn}
    \mu_{\tk,N}
    := \frac{1}{N}\sum_{\ell=1}^{N}\delta[x_{\tk,\ell}]
    \approx \mu_{\tk}.
\end{equation}
We write $\mu_{\tk,N} \in P_{2,N}(\mathbb{R}^{d})$, where
\begin{equation*}
    P_{2,N}(\mathbb{R}^{d})
    :=
    \bigg\{
        \frac{1}{N}\sum_{\ell=1}^{N}\delta[x_{\ell}]
        : x_{\ell}\in\mathbb{R}^{d}
    \bigg\}
    \subset P_{2}(\mathbb{R}^{d})
\end{equation*}
denotes the set of equally weighted atomic measures with $N$ particles. In one-dimension, we additionally require an order of particles $x_{\ell} \le x_{\ell + 1}$ for all $\ell$.

This leads to the following spatially discrete JKO scheme.
\begin{definition}[Particle method]\rm
Given $\mu_{\tk,N}\in P_{2,N}(\mathbb{R}^{d})$, we define $\mu_{\tk+1,N}$ by
\begin{equation} \label{1001}
    \mu_{\tk+1,N}
    \in
    \operatorname*{argmin}_{\mu_N\in P_{2,N}(\mathbb{R}^{d})}
    \left\{
        \mathcal{F}_{N}(\mu_N)
        + \frac{1}{2\tau}
        W_{2}^{2}(\mu_N,\mu_{\tk,N})
    \right\},
\end{equation}
where $\mathcal{F}_{N}$ denotes a suitable particle approximation of the free energy $\mathcal{F}$.
\end{definition}

Recall from~\eqref{2} that the free energy consists of the internal,
confinement, and interaction energies. For continuous confinement
and interaction potentials $V$ and $W$, the corresponding
energies can be evaluated directly on an atomic measure $\mu_N=\frac{1}{N}\sum_{\ell=1}^{N}\delta[x_\ell]$:
\begin{equation} \label{F_VN}
    \mathcal{F}_{V}(\mu_N)
    =
    \frac{1}{N}\sum_{\ell=1}^{N}V(x_\ell), 
    \qquad 
    \mathcal{F}_{W}(\mu_N)   =  \frac{1}{2N^{2}}  \sum_{\ell,m=1}^{N} W(x_\ell-x_m).
\end{equation} 
The \textit{principal challenge}, however, lies in the internal energy $\mathcal{F}_{U}(\mu)$, which is generally not well defined on singular atomic measures, and consequently, a suitable particle approximation of $\mathcal{F}_{U}$ is required. In the following, we consider three different constructions. The first exploits the ordering structure available specifically in one spatial dimension, while the other two apply more generally in arbitrary dimensions.

\subsubsection{One-dimensional}
On $\mathbb{R}$, the structure of transport map can be explicitly decomposed via the increasing rearrangement of mass (see, for example, \citep[Chapter 1]{villani2008optimal}). For any measure $\mu \in P_{2}(\mathbb{R})$, let $F_{\mu}(t) = \mu\big((-\infty, t)\big)$ be its distribution function, we define the quantile map: 
\begin{equation*}
    X_{\mu}: [0, 1] \to \mathbb{R} \cup \{\pm \infty\}, \qquad X_{\mu}(\theta) := \sup\Big\{t \in \mathbb{R} : F_{\mu}(t) \le \theta\Big\}.
\end{equation*}
Then $X_\mu$ is the left-continuous pseudo inverse of $F_{\mu}(t)$ and $(X_\mu)_\# \nu_{[0,1]} = \mu$, where $\nu_{[0,1]}$ is the Lebesgue measure over $[0,1]$. Likewise, let $X_k$ be the quantile of $\mu_k$, and then $(X_k)_\# \nu_{[0,1]} = \mu_k$. Therefore \eqref{1} can be rewritten as 
\begin{equation*}
    \min_{X} \mathcal{F}(X_\# \nu_{[0, 1]}) + \frac{1}{2\tau} \int_{0}^1 |X(\theta) - X_{\tk}(\theta)|^{2} \ud \nu_{[0, 1]}(\theta),
\end{equation*}
where $X$ is in the space of non-decreasing, $L^{2}$ functions on $[0, 1]$.

In practice, we do not work directly with $X$, but consider its cell average:
\begin{equation*}
    x_{\ell} = N\int_{(\ell - 1)/N}^{\ell/N} X(\theta) \ud \theta,\quad \text{ for $\ell = 1, \cdots, N$}\,.
\end{equation*}
We then define the empirical measure $\mu_N:=\frac1N\sum_{\ell=1}^N\delta[x_\ell]$.
Writing $U_{\rm mc}(s):=s\,U(1/s)$,
the internal energy can be expressed in Lagrangian coordinates as
\[
\mathcal F_U(\mu)=\int_0^1 U_{\rm mc}\bigl(X'(\theta)\bigr)\ud\theta.
\]
This identity is understood with the recession-function convention $U_{\rm mc}(0)=U_\infty$. In particular, $U_\infty=+\infty$ for $m\ge1$, whereas $U_\infty=0$ for $m<1.$
To obtain a discrete approximation, we consider a piecewise-affine interpolation of particle locations $x_\ell$. Since neighboring particles correspond to reference points separated by $1/N$, the slope of this interpolation between $x_\ell$ and $x_{\ell+1}$ is
\[
X_N'(\theta)
=
\frac{x_{\ell+1}-x_\ell}{1/N}
=
N(x_{\ell+1}-x_\ell).
\]
Therefore, discretizing the Lagrangian internal energy gives 
\[
\mathcal F_{U, N}(\mu_N):=\int_0^1 U_{\rm mc}\bigl(X_N'(\theta)\bigr)\ud\theta  \approx \sum_{\ell=0}^{N-1} \frac1N U_{\rm mc}\bigl(N(x_{\ell+1}-x_\ell)\bigr) = \sum_{\ell=0}^{N-1} U\left(\frac1{N(x_{\ell+1}-x_\ell)}\right)(x_{\ell+1}-x_\ell).
\]
This type of discretization has been used, for instance, in \cite{di2015rigorous}. A minor subtlety arises when two or more particles coincide, that is, when $x_\ell=x_{\ell+1}$ for some $\ell$, since the expression above is then not directly well defined. To accommodate this case, we adopt a slightly modified definition based on the value of
$\lim_{t\to\infty}\frac{U(t)}{t}$,
as specified in~\eqref{3-9-internal} below.

The quantile formulation has motivated several variational Lagrangian and deterministic particle methods in one dimension. For example, Matthes and Osberger \cite{matthes2014convergence} proved convergence of a fully discrete variational scheme for nonlinear drift--diffusion on bounded intervals as the spatial mesh and time step vanish under an inverse CFL condition. Daneri, Radici, and Runa \cite{daneri2022deterministic} established convergence of continuous-time deterministic particle approximations for aggregation--diffusion on the real line, while Carrillo, Patacchini, Sternberg, and Wolansky \cite{carrillo2016convergence} proved convergence of diffusive particle gradient flows on bounded intervals using a different internal-energy reconstruction. At the level of time-discrete variational schemes, Di Marino, Portinale, and Radici \cite{di2024optimal} established fixed-time-step convergence for particle approximations of transport distances with nonlinear mobilities. Their JKO convergence theorem assumes continuous convergence of the discrete energies in $p$-Wasserstein metrics with $p < 2$, a condition that does not generally hold for particle approximations of nonlinear internal energies.

Our focus is the spatial approximation of the JKO minimization problem itself, with the time step held fixed. We study whether the adjacent-gap discretization of the internal energy, combined with the atomic Wasserstein transport cost, yields consistent particle minimizers that converge to the continuum JKO minimizer on $\mathbb{R}$. This requires more than convergence of empirical measures: the nonlinear internal energy must be handled through a variational approximation, and the second-moment bounds supplied by coercivity yield only narrow compactness, not compactness in $W_2$ topology. We address these issues by establishing $\Gamma$-convergence of the discrete JKO functionals and upgrading convergence of their minimizers to $W_2$-convergence. Under the standing assumptions, including uniqueness of the limiting JKO minimizers, this convergence can then be propagated through successive time steps. The analysis accommodates both internal energies and confinement and interaction potentials whose convex components may grow superquadratically, and separates spatial consistency from the vanishing-time-step limit. The details are presented in Section~\ref{sec:3}.

\subsubsection{Blob method}
Another way to handle the singularity of atomic measures is to regularize them by convolution with a mollifier, leading to the blob method \cite{craig2016blob, craig2016convergence, liu2017random, carrillo2019aggregation, carrillo2019blob}. More specifically, one replaces the atomic measure $\mu_N$ by the continuous density $\psi_\varepsilon * \mu_N$ and approximates the internal energy by the well-defined functional
\begin{equation}\label{2-bob}
    \mathcal{F}_{\varepsilon,U}(\mu_N)
    :=
    \int_{\mathbb{R}^{d}}
    U(\psi_{\varepsilon} \ast \mu_N)\,\ud x.
\end{equation}
Here $  \psi_{\varepsilon}(x) :=
    \varepsilon^{-d}\psi(x/\varepsilon)$,
where $\psi$ is a continuous and bounded probability density with finite second moment. The convolution of $\psi_\varepsilon$ with a measure $\mu$ is defined by
\begin{equation*}
    (\psi_{\varepsilon} \ast \mu)(x)
    :=
    \int_{\mathbb{R}^{d}}
    \psi_{\varepsilon}(x-y)\,\ud\mu(y).
\end{equation*}
Correspondingly, we consider the regularized version of \eqref{1} for blob method: 
\begin{align} \label{JKO_eps}
    \JKO_\varepsilon (\mu,\mu_{\tk}):=   \mathcal{F}_{\varepsilon, U}(\mu) + \mathcal F_V(\mu) + \mathcal F_W(\mu)
        + \frac{1}{2\tau} W_{2}(\mu,\mu_{\tk})^{2}\,.
\end{align}

This energy regularization was introduced by Carrillo, Craig, and Patacchini in \cite{carrillo2019blob}. They established $\Gamma$-convergence of the regularized internal energies for $m\ge1$ and, for $m\ge2$, convergence of the associated continuous-time Wasserstein gradient flows to solutions of \eqref{0} as $\varepsilon\to0$, under suitable assumptions on the potentials, initial data, and regularized solutions. Their particle approximation is formulated as a system of ordinary differential equations, and its convergence is established at the spatially discrete, continuous-in-time level, with initial particle spacing $h=o(\varepsilon)$. 

In contrast, our implicit particle approach starts from the time-discrete JKO formulation~\eqref{JKO_eps} and approximates the resulting variational minimization problem using particles. In Section~\ref{sec:4}, we establish the convergence of this particle approximation based on a specific choice of mollifier, with the regularization parameter $\varepsilon>0$ held fixed. We also establish simultaneous convergence as the particle number tends to infinity and the regularization parameter vanishes, under the scaling $N\varepsilon_N^d\to\infty.$

\subsubsection{Transport map representation}
A third approach to handling~\eqref{1001} is based on a transport-map representation. By Brenier's theorem, under suitable assumptions on $\mu_{\tk}$, the minimization problem~\eqref{1} can be reformulated in terms of a convex potential $\phi$ as
\begin{equation}\label{def:jko-b}
    \min_{\substack{\phi:\mathbb{R}^{d}\to\mathbb{R}\, \mathrm{convex}}}
    \JKO\big((\nabla\phi)_{\#}\mu_{\tk},\mu_{\tk}\big)
    =
    \mathcal{F}\big((\nabla\phi)_{\#}\mu_{\tk}\big)
    +
    \frac{1}{2\tau}
    \int_{\mathbb{R}^{d}}
    |\nabla\phi(x)-x|^{2}\,\ud\mu_{\tk}(x).
\end{equation}
Motivated by this representation, the corresponding particle formulation is
\begin{equation}\label{def:jko-bN}
    \min_{\substack{\phi:\mathbb{R}^{d}\to\mathbb{R} \;\mathrm{convex}}}
    \JKO_N\big(\phi;\mu_{\tk,N}\big)
    :=
    \mathcal{F}_{N}\big(\phi;\mu_{\tk,N}\big)
    +
    \frac{1}{2\tau}
    \int_{\mathbb{R}^{d}}
    |\nabla\phi(x)-x|^{2}\,\ud\mu_{\tk,N}(x).
\end{equation}
In practice, the convex potential $\phi$ can be parameterized by a neural network, and such transport-map-based approaches have attracted considerable recent attention
\cite{makkuva2020optimal, korotin2022neural, lee2024deep, feng2026learn, dumont2026learning}.

It remains to specify how the internal energy $\mathcal{F}_{U}$ is evaluated in~\eqref{def:jko-bN}. We follow the approach in \cite{lee2024deep}. To motivate the construction, we let $\nu$ be an absolutely continuous measure and suppose that $\nabla\phi$ is sufficiently regular and invertible. By the change-of-variables formula,
\begin{equation*}
    \mathcal{F}_{U}\big((\nabla\phi)_{\#}\nu\big)
    =
    \int_{\mathbb{R}^{d}}
    S\left(
        \frac{\nu(x)}
        {\det(\nabla^{2}\phi(x))}
    \right)
    \ud\nu(x),
    \qquad
    S(s):=\frac{U(s)}{s}.
\end{equation*}
If $\nu$ is replaced directly by an atomic measure, integration with respect to $\ud\nu$ remains meaningful, but the density value $\nu(x)$ appearing inside $S$ is no longer defined. The key idea is therefore to use the particle measure for quadrature while evaluating the density through an associated absolutely continuous measure.

More precisely, for $\mu_{\tk, N}$ of the form \eqref{mu_kn}, 
we define the discrete internal energy by
\begin{equation*}
    \mathcal{F}_{U,N}
    \big(
        \phi;
        \mu_{\tk,N}^{\mathrm{ac}},
        \mu_{\tk,N}
    \big)
    :=
    \frac{1}{N}
    \sum_{\ell=1}^{N}
    S\left(
        \frac{
            \mu_{\tk,N}^{\mathrm{ac}}(x_{\tk,\ell})
        }{
            \det\big(\nabla^{2}\phi(x_{\tk,\ell})\big)
        }
    \right).
\end{equation*}
Here, $\mu_{\tk,N}^{\mathrm{ac}}$ is obtained by transporting the initial absolutely continuous measure $\mu_{0}$ through the sequence of maps computed at the preceding time steps. More precisely, if $\phi_{\tk,N}$ denotes the transport potential obtained at the $\tk$-th step, then
\begin{equation*}
    \mu_{\tk,N}^{\mathrm{ac}}
    :=
    \big(
        \nabla\phi_{\tk,N}
        \circ
        \nabla\phi_{\tk-1,N}
        \circ\cdots\circ
        \nabla\phi_{1,N}
    \big)_{\#}\mu_{0}.
\end{equation*}
At the same time, the particle locations are transported by the same maps according to
\begin{equation*}
    x_{\tk,\ell}
    =
    \nabla\phi_{\tk,N}
    (x_{\tk-1,\ell}),
    \qquad
    \ell=1,\ldots,N,
\end{equation*}
which yields the empirical measure $ \mu_{\tk,N}=
    \frac{1}{N}
    \sum_{\ell=1}^{N}
    \delta[x_{\tk,\ell}]$.
Thus, the absolutely continuous measure $\mu_{\tk,N}^{\mathrm{ac}}$ and the empirical measure $\mu_{\tk,N}$ evolve under the same sequence of transport maps: the former supplies the density information needed to evaluate the internal energy, while the latter supplies the density values needed to evaluate the internal energy through the Jacobian determinant of the transport map. The convergence analysis therefore needs to account simultaneously for the approximation of the source measure, the variation of the optimizing transport potential, and the behavior of the energy integrands at spatial infinity, see Section~\ref{sec:5} below. With the time step fixed, we establish variational convergence of the particle-discretized JKO problem as the number of particles tends to infinity, together with $W_2$-convergence of the resulting particle measures. The unbounded-domain setting requires additional tail estimates, which are not needed on bounded domains and are essential for propagating the convergence through successive JKO steps.

We also mention that ransport-map formulations have also motivated variational discretizations in multiple dimensions on bounded domains. For example, Benamou, Carlier, M\'erigot, and Oudet \cite{benamou2016discretization} considered two bounded convex domains $X,Y\subset\mathbb{R}^d$, with the source measure supported on $X$ and the transported measure constrained to $Y$. They constructed a geometric discretization based on the Monge--Amp\`ere operator and established convergence of the corresponding minimization problems. Junge, Matthes, and Osberger \cite{JungeMatthesOsberger2017} considered nonlinear Fokker--Planck equations on the multidimensional unit cube $\Omega=(0,1)^d$ and developed a fully discrete variational scheme in which the incremental transport maps are approximated in finite-dimensional spaces of gradients of smooth basis functions. They proved convergence to classical solutions as both the time step and the spatial resolution of the transport-map approximation vanish, under uniform regularity and positivity assumptions on the discrete densities. In their numerical implementation, the integrals in the variational problem are further approximated by quadrature, although this additional quadrature approximation is not included in their convergence analysis.

\subsection{Main result}
\begin{assumption}\label{assumption}\rm Throughout the paper, we assume that $U, V, W$ in~\eqref{0} satisfies the following.
\begin{enumerate}
    \item For $m > d/(2 + d)=:m_0$, the internal energy takes the form
    \begin{equation*}
        U_{m}(s) =
        \begin{cases}
            s\log(s) & m = 1\\
            \frac{1}{m - 1}s^{m} & m \neq 1
        \end{cases}
    \end{equation*}
    We often omit the subscript $m$ in $U_{m}$, and denote it to be $U$ to simply the notation. 
    
    \item $V, W: \mathbb{R}^{d} \to \mathbb{R}$ are continuous functions that can be decomposed as follows,
    \begin{equation*}
        V = V_{\mathrm{pos}} + V_{\mathrm{res}}, \quad W = W_{\mathrm{pos}} + W_{\mathrm{res}}
    \end{equation*}
    where $V_{\mathrm{pos}}$ and $W_{\mathrm{pos}}$ are continuous, positive and convex, and $V_{\mathrm{res}}$ and $W_{\mathrm{res}}$ are sub-quadratic in the sense that  for some positive constant $C$ and $0 \le \beta < 2$
    \begin{equation*}
        |V_{\mathrm{res}}(x)| \le C\big(1 + |x|^{\beta}\big),\quad |W_{\mathrm{res}}(x)| \le C\big(1 + |x|^{\beta}\big)
    \end{equation*}
\end{enumerate}
\end{assumption}

For the internal energy, the constant $m_{0}$ is the critical exponent ensuring that $\mathcal{F}_{U}$ is well defined on $P_{2}(\mathbb{R}^{d})$; see Lemma~\ref{lemma-5-internal} below. For the confinement and interaction potentials $V$ and $W$, we allow each potential to be decomposed into two parts: a convex part and a subquadratic part. The decomposition enables us to accommodate a broad class of potentials with different growth and convexity properties.
The subquadratic components give rise to energies that are continuous with respect to convergence in $P_{2}(\mathbb{R}^{d})$, while still allowing for certain nonconvex potentials. The convex components, on the other hand, allow for potentials with superquadratic growth, such as $ V(x)=|x|^{p},
~
    W(x)=|x|^{p}$
with arbitrarily large $p$. In particular, when $p>2$, the corresponding potential energies are generally not continuous with respect to convergence in $P_{2}(\mathbb{R}^{d})$, and therefore require a different treatment from their subquadratic counterparts.

\begin{assumption}[uniqueness]\rm We assume the solutions of the spatially continuous JKO scheme~\eqref{1}, and the blob-regularized JKO scheme \eqref{JKO_eps} are unique. 
\end{assumption}

The existence of solutions to the JKO scheme follows directly from the stability result in Lemma~\ref{lemma-7-stability}, together with the moment lower bounds established in Lemmas~\ref{lemma-5-internal} and~\ref{lemma-6-well} below. For the internal energy $U_m$, McCann's displacement--convexity condition holds when $m\ge 1-1/d$. Consequently, under this additional restriction, if $V$ and $W$ are convex, the free energy $\mathcal F$ is displacement convex. The uniqueness of JKO minimizers can then be obtained under the corresponding convexity hypotheses. For $d/(d+2)<m<1-1/d$, which is possible when $d>2$, uniqueness remains part of Assumption 2 \cite[Section~3]{carrillo2003kinetic} \cite[Chapter~4]{ambrosio2008gradient}.

Alternatively, for $m\geq 1$ and positive potentials $V, W$, uniqueness can be established through convexity with respect to linear interpolation of measures. In this case, $\mathcal{F}_{U}$ is strictly convex, while $W_2^2(\cdot,\nu)$ and $\mathcal{F}_{V}$ are convex under linear interpolation. If, in addition, the interaction energy $\mathcal{F}_{W}$ is convex under linear interpolation under the corresponding assumptions on $W$, then the JKO functional is strictly convex and therefore admits at most one minimizer; see \cite{JKO}.

\begin{assumption} \label{assump3}
For the \textit{transport-map approach}, we only consider subquadratic potentials $V$ and $W$, and assume that the initial density $\mu_0\in C(\mathbb{R}^d)\cap P_2(\mathbb{R}^d)$ and satisfies $\|\mu_0\|_{L^{\infty}}\leq M$ with some constant $M>0$. In addition, depending on the value of $m$, we impose the following lower bounds on the initial condition:
\begin{equation}\label{lower-bound-3}
    \begin{cases}
        \displaystyle
        \mu_0(x)\geq \kappa
        \exp\!\left(-\frac{|x|^2}{4\gamma}\right),  
        & m=1,\\[2mm]
        \mu_0(x)>0,
        & m>1,\\[2mm]
        \displaystyle
        \mu_0(x)\geq
        \big(Z+\kappa |x|^2\big)^{-1/(1-m)},
        & 1-\dfrac{2}{2+d}<m<1,
    \end{cases}
\end{equation}
for some positive constants $\kappa$, $\gamma$, and $Z$.

We further suppose that, at each JKO step $\tk$, the Brenier potential associated with the minimizer of~\eqref{def:jko-b} belongs to $ C_{\Lambda_{\tk}}^{2,\alpha_{\tk}}(\mathbb{R}^d)$
for some constants $\Lambda_{\tk}, \alpha_{\tk} \in \big(0, )$. Here, for $\Lambda, \alpha \in (0,1)$, we define
\begin{equation}\label{Calpha}
    C_{\Lambda}^{2, \alpha}(\mathbb{R}^{d}) := \Big\{\phi \in C^{2}(\mathbb{R}^{d}): \Lambda I_{d} \le \nabla^{2}\phi(x) \le \tfrac{1}{\Lambda}I_{d},~ [\nabla^{2}\phi]_{\alpha; B_{R}(0)} \le C_{R} \text{ for all $R > 0$}\Big\},
\end{equation}
where $C_R$ may depend on $R$, and
\begin{equation*}
    [f]_{\alpha;B_R(0)}
    :=
    \sup_{\substack{x,y\in B_R(0)\\x\neq y}}
    \frac{\lVert f(x)-f(y)\rVert}{|x-y|^\alpha}.
\end{equation*}
\end{assumption}


For degenerate diffusion with $m>1$, compactly supported initial data $\mu_0$ are generally admissible for~\eqref{0}, and the corresponding solution $\mu(t)$ may remain compactly supported due to the finite speed of propagation characteristic of degenerate diffusion. In the present work, however, we impose the stronger assumption that $\mu_0$ is strictly positive on $\mathbb{R}^d$. This condition is needed to ensure uniqueness of the Brenier potential within the class of convex functions on $\mathbb{R}^d$. Indeed, if $\mu_{\tk}$ is compactly supported, the associated Brenier map is determined only $\mu_{\tk}$-almost everywhere, and the Brenier potential is therefore generally nonunique outside $\operatorname{supp}(\mu_{\tk})$. By the standard maximum principle, strict positivity is preserved by~\eqref{0}, so that the solution of the flow $\mu(t)$ remains strictly positive for positive times whenever $\mu_0$ is strictly positive.

In the case of linear diffusion $(m = 1)$, the $C_{\Lambda}^{2, \alpha}(\mathbb{R}^{d})$ assumption of the Brenier potential can be guaranteed for initial value with Gaussian decay. In particular, one can prove the persistence of tail, which is, if we started from an initial value with Gaussian decay, the solution remains Gaussian decay for all time, can be proved via the composition principal of parabolic PDEs based on a barrier construction as that in \cite{kim2012patlak}. This persistence of tail condition ensures that the Brenier potential $\phi$ grows like $|x|^{2}$ at infinity, and implies that $\nabla^{2}\phi$ remains uniformly elliptic at infinity. Based on the recent regularity results of Coudreuse \cite{coudreuse2025li, coudreuse2026aronson}, one can prove that the Brenier potential $\phi_{\tk}$ falls into the class $C_{\Lambda}^{2, \alpha}(\mathbb{R}^{d})$ with some $\alpha$ and $\Lambda$.

\begin{definition}[Recovery sequence] \label{defn:recoveryS} Let $\mathcal F_N,\mathcal F: P_2(\mathbb R^d)\to(-\infty,+\infty]$. A sequence $(\mu_N)_N \subset P_2(\mathbb R^d)$ is a recovery sequence for $\mu\in P_2(\mathbb R^d)$ if it converges narrowly to $\mu$ as $N$ tends to $+\infty,$ and 
\[
\limsup_{N \to +\infty} \mathcal F_N(\mu_N)\le\mathcal F(\mu).
\]
It is a $W_2$-convergent recovery sequence if, in addition, $\big(W_2(\mu_N,\mu)\big)_N$ converges to $0.$  
\end{definition}

We now summarize our main results. 

\begin{theorem}[1D]\label{theorem:real}Given initial value $\mu_{0} \in P_{2}(\mathbb{R})$ and let $\mu_{0, N} := 1/N\sum_{\ell = 1}^{N} \delta[x_{0, \ell}] \in P_{2, N}(\mathbb{R})$ be such that $(\mu_{0, N})_N$ converges to $\mu_{0}$ in $P_{2}(\mathbb{R}^{d})$. Then $(\mu_{\tk, N})_N$ converges to $\mu_{\tk}$ in $P_{2}(\mathbb{R}^{d}),$ for all $\tk \in \mathbb{N}.$ Here, $\mu_{\tk, N}$ is the $N$-particle solution \eqref{pJKO-1D}, and $\mu_{\tk}$ is the solution of JKO scheme \eqref{1}.
\end{theorem}
\begin{theorem}[blob method]\label{theorem:blob:bis} Suppose that $(\mu_{0,N})_N\subset  P_{2,N}(\mathbb R^d)$ converges to $\mu_0$ in the $W_2$ metric,  and let 
$\mu_{k,N}$ be generated recursively by the particle minimization problem with $\varepsilon=\varepsilon_N$. Suppose that $m>d/(d+2),$ $ \psi\in C_c^\infty(\mathbb R^d)$ is a probability density, suppose $\big(\varepsilon_N\big)_N$ decreases to $0$ and $\big(N\varepsilon_N^d\big)_N$ tends to $+\infty$ as $N$ tends to $+\infty.$ Define
\[
\mathcal F_N(\mu)=
\begin{cases}
\mathcal F_{U}((\psi_{\varepsilon_N}*\mu) \ud x)+\mathcal F_V(\mu)+\mathcal F_W(\mu),&\mu\in P_{2,N}(\mathbb R^d),\\
+\infty,&\text{otherwise}.
\end{cases}
\] 
Then, for every sequence $\mu_N\rightharpoonup\mu$ narrowly with uniformly bounded second moments,
\[
\liminf_{N \to +\infty}\mathcal F_N(\mu_N)\ge\mathcal F(\mu).
\] 
Moreover, every $\mu\in P_2(\mathbb R^d)$ admits a $W_2$-convergent recovery sequence. In this sense, $$\Gamma\text{-}\lim_{N \to \infty}\mathcal F_N = \mathcal F.$$
Furthermore,  under the uniqueness of the limiting JKO minimizers, for fixed $\tau>0$ and $\tk$, $(\mu_{\tk,N})_N$ converges to $\mu_\tk$ in $ P_{2,N}(\mathbb R^d)$.
Here $\mu_{\tk,N}$ solves the implicit particle problem with $\varepsilon=\varepsilon_N$, whereas $\mu_\tk$ solves the unregularized JKO problem.
\end{theorem} 

\begin{remark}[Conjectured optimality of the scaling condition] Let $m>1$, and $(\varepsilon_N)_N \subset (0,+\infty)$ decreasing to $0$.  Let $\mathcal F_N$ be as in Theorem \ref{theorem:blob:bis}. We conjecture that 
\[
\Gamma\text{-}\lim_{N \to \infty}\mathcal F_N = \mathcal F \quad\Longleftrightarrow\quad N\varepsilon_N^d\to+\infty,
\]
where $\Gamma$-convergence is understood  in the sense of Theorem \ref{theorem:blob:bis}, including the existence of $W_2$--convergent recovery sequences. Thus, we conjecture that the scaling condition in Theorem \ref{theorem:blob:bis} is optimal here. 
\end{remark}

\begin{theorem}[transport map method]\label{theorem-3-neural}Given initial value $\mu_{0} \in P_{2}(\mathbb{R}^d)$ and let $\mu_{0, N} := \frac{1}{N}\sum_{\ell = 1}^{N} \delta[x_{0, \ell}] \in P_{2, N}(\mathbb{R}^d)$ be a sampling of $\mu_{0}$ such that $\mu_{0, N} \to \mu_{0}$ in $P_{2}(\mathbb{R}^{d})$. Then there exists a sequence of minimizers $(\phi_{\tk, N})_{\tk = 1}^{\infty}$ to the discrete JKO loss functional \eqref{1002} under the constraint $C_{\Lambda}^{2, \alpha}(\mathbb{R}^{d})$, and we have $\nabla \phi_{\tk, N} \to \nabla \phi_{\tk}$ locally uniformly as $N \to \infty$ for all $\tk$ where $\phi_{\tk}$ is Brenier potential solution at the $\tk$-th step, and
\begin{equation*}
    \mu_{\tk, N} \to \mu_{\tk} \text{ in } P_{2}(\mathbb{R}^{d}) \text{ as } N \to \infty \text{ for all } \tk \in \mathbb{N}
\end{equation*}
where $(\mu_{\tk})_{\tk = 1}^{\infty}$ is the solution of the JKO scheme \eqref{1}. 
\end{theorem}

\section{Preliminary}\label{Sec:2}

\subsection{Optimal transport}
In this part, we provide some preliminary discussions on the optimal transport theory. For a thorough investigation, we refer readers to \cite{ambrosio2008gradient, villani2008optimal, santambrogio2015optimal}.

\subsubsection{Basic concepts}
Throughout this paper, we denote by $\nu_{[0,1]^d}$ the Lebesgue measure restricted to $[0,1]^d$. We write $P(\mathbb{R}^d)$ for the space of probability measures on $\mathbb{R}^d$, and $P^{\mathrm{ac}}(\mathbb{R}^d)$ for the subset of probability measures that are absolutely continuous with respect to the Lebesgue measure.

Given a measurable map $\br:\mathbb{R}^d\to\mathbb{R}^d$ and a measure $\mu\in P(\mathbb{R}^d)$, the push-forward of $\mu$ under $\br$, denoted by $\br_{\#}\mu$, is defined by $(\br_{\#}\mu)(A):= \mu\big(\br^{-1}(A)\big)$
for every Borel set $A\subset\mathbb{R}^d$. We say that $\br$ transports $\mu$ to $\nu$ if $\nu=\br_{\#}\mu$.

For $p\geq1$, we define respectively  the $p$-th moment of $\mu \in P^{\mathrm{ac}}(\mathbb{R}^d)$ and the $p$-Wasserstein space by  
\begin{equation*}
M_p(\mu):=\int_{\mathbb{R}^d}|x|^p\,\ud\mu(x), \quad 
    P_p(\mathbb{R}^d)
    :=
    \left\{
        \mu\in P(\mathbb{R}^d): M_p(\mu)<\infty
    \right\}.
\end{equation*}
The space $P_p(\mathbb{R}^d)$ is equipped with a metric called the $p$-Wasserstein distance. In particular, for $p=2$,
\begin{equation*}
    W_2^2(\mu,\nu)
    :=
    \min_{\gamma\in\Pi(\mu,\nu)}
    \int_{\mathbb{R}^d\times\mathbb{R}^d}
    |x-y|^2\,\ud\gamma(x,y),
\end{equation*}
where
\begin{equation*}
    \Pi(\mu,\nu)
    :=
    \left\{
        \gamma\in P(\mathbb{R}^d\times\mathbb{R}^d):
        (\operatorname{proj}_x)_{\#}\gamma=\mu,\ 
        (\operatorname{proj}_y)_{\#}\gamma=\nu
    \right\}
\end{equation*}
is the set of couplings between $\mu$ and $\nu$. We denote an optimal coupling by $\gamma_{\mathrm{o}}\in\Pi(\mu,\nu)$. This coupling-based definition is known as the Kantorovich formulation of optimal transport and, in particular, remains well defined when both $\mu$ and $\nu$ are atomic.

When the source measure $\mu$ is absolutely continuous, the Kantorovich formulation is equivalent to the Monge formulation through Brenier's theorem \cite{brenier1987decomposition, GangboMcCann1996}, as recalled below.

\begin{lemma}[Brenier's theorem]
Suppose that $\mu\in P_2^{\mathrm{ac}}(\mathbb{R}^d)$ and $\nu\in P_2(\mathbb{R}^d)$. Then there exists a convex function
$\phi_{\mu\to\nu}:\mathbb{R}^d\to\mathbb{R}$ such that
$\nabla\phi_{\mu\to\nu}$ is the optimal transport map from $\mu$ to $\nu$. In particular,
\begin{equation*}
    \nu   =    (\nabla\phi_{\mu\to\nu})_{\#}\mu,    \qquad    W_2^2(\mu,\nu)    =    \int_{\mathbb{R}^d}    |\nabla\phi_{\mu\to\nu}(x)-x|^2\,\ud\mu(x),
\end{equation*}
and the corresponding optimal coupling is $\gamma_{\mathrm{o}} = (\mathrm{id},\nabla\phi_{\mu\to\nu})_{\#}\mu.$ 
Moreover, the optimal transport map $\nabla\phi_{\mu\to\nu}$ is unique $\mu$-almost everywhere.
\end{lemma}

\subsubsection{Convergence in the space of probability measures}
In the following, we collect several facts concerning the topology of spaces of probability measures.

\begin{definition}[Narrow convergence]\rm
Let $(\mu_N)_{N=1}^{\infty}$ be a sequence of finite positive Radon measures on $\mathbb{R}^d$. We say that $(\mu_N)_N$ converges \emph{narrowly} to a measure $\mu$, and write
\[
    \mu_N \rightharpoonup \mu
    \qquad\text{narrowly as }N\to\infty,
\]
if
$\lim_{N\to\infty}
    \int_{\mathbb{R}^d}\varphi\,\ud\mu_N
    =
    \int_{\mathbb{R}^d}\varphi\,\ud\mu$ for every 
$\varphi\in C_b(\mathbb{R}^d)$,
where $C_b(\mathbb{R}^d)$ denotes the space of bounded continuous functions on $\mathbb{R}^d$.
\end{definition}

\begin{definition}[Tightness]\rm
A sequence of finite positive Radon measures $(\mu_N)_{N=1}^{\infty}$ on $\mathbb{R}^d$ is said to be \emph{tight} if, for every $\varepsilon>0$, there exists $R>0$ such that $ \mu_N\big(\mathbb{R}^d\setminus B_R(0)\big)<\varepsilon$ for all $N$.
\end{definition}

The following Prokhorov's theorem provides a compactness criterion for probability measures.
\begin{lemma}[Prokhorov's theorem]
Let $(\mu_N)_{N=1}^{\infty}\subset P(\mathbb{R}^d)$ be a tight sequence of probability measures. Then there exist a subsequence, not relabeled, and a measure $\mu\in P(\mathbb{R}^d)$ such that $(\mu_N)_N$ converges \emph{narrowly} to a measure $\mu.$ 
In particular, for a sequence $(\mu_N)_{N=1}^{\infty}\subset P_2(\mathbb{R}^d)$, tightness can be guaranteed by a uniform bound on the second moments: $\sup_N M_2(\mu_N)<\infty$.
\end{lemma}
%
%
\begin{definition}[$W_2$-convergence]\rm
Let $(\mu_N)_{N=1}^{\infty}\subset P_2(\mathbb{R}^d)$ and $\mu\in P_2(\mathbb{R}^d)$. We say that $(\mu_N)_N$ converges to $\mu$ in $P_2(\mathbb{R}^d)$, or equivalently in $W_2$, if $\big(W_2(\mu_N,\mu)\big)_N$ converges to $0.$
\end{definition}

The topology induced by $W_2$ is stronger than the narrow topology: convergence in $W_2$ is equivalent to narrow convergence together with convergence of the second moments. More precisely, we have the following equivalent characterizations; see \cite[Definition~6.8]{villani2008optimal}.

\begin{lemma}\label{lemma-3-p2}
Let $(\mu_N)_{N=1}^{\infty}\subset P_2(\mathbb{R}^d)$ and $\mu\in P_2(\mathbb{R}^d)$. Then, $(\mu_N)_N$ converges to $\mu$ in $P_2(\mathbb{R}^d)$ 
if and only if $\mu_N\rightharpoonup\mu$ narrowly and any one of the following equivalent conditions holds:
\begin{enumerate}
    \item $\lim_{N\to\infty}M_2(\mu_N)=M_2(\mu)$.

    \item The second moments are uniformly integrable, namely, for every $\varepsilon>0$, there exists $R>0$ such that $\int_{|x|\geq R}|x|^2\,\ud\mu_N(x)<\varepsilon$ for all $N$.

    \item For every continuous function $\varphi:\mathbb{R}^d\to\mathbb{R}$ satisfying
$ |\varphi(x)|\leq C(1+|x|^2)$
    for some constant $C>0$, one has $ \int_{\mathbb{R}^d}\varphi\,\ud\mu_N
        \to
        \int_{\mathbb{R}^d}\varphi\,\ud\mu$.
\end{enumerate}
\end{lemma}

\subsubsection{Canonical sampling}
Finally, we provide a canonical way of constructing $N$-particle approximation of any $\mu \in P_{2}(\mathbb{R}^{d})$. Let $(E_{\ell})_{\ell = 1}^{N}$ be an $N$-partition of the unit cube $[0, 1]^{d}$ given $N \in \mathbb{N}$, such that $E_{\ell}$ are connected, and
\begin{equation}\label{2-1-3-1}
\lim_{N \to +\infty}   \Big( \max_{1 \le \ell \le N} \mathrm{diam}(E_{\ell})\Big) =0,\qquad \nu_{[0, 1]^{d}}(E_{\ell}) = N^{-1}.
\end{equation}
For real-line $\mathbb{R}$, we partition the reference domain $[0, 1]$ into $N$ uniform cells $[(\ell - 1)/N, \ell/N]$.

\begin{definition}[canonical $N$-particle approximation]\rm\label{def-7}
Let $\mu \in P_{2}(\mathbb{R}^{d})$ and $\phi_{\to \mu}$ be the Brenier potential from $\nu_{[0, 1]^{d}}$ to $\mu$. We define the corresponding $N$-particle approximation as follows,
\begin{equation} \label{eqn:xl}
    \mu_{N} := \frac{1}{N} \sum_{\ell = 1}^{N} \delta[x_{\ell}] \quad \text{ where } \quad x_{\ell} = \frac{1}{\nu_{[0, 1]^{d}}(E_{\ell})}\int_{E_{\ell}} \nabla \phi_{\to \mu}(\theta) \ud \nu_{[0, 1]^{d}}.
\end{equation}
\end{definition}

\begin{lemma}\label{lemma-7-4}For any $\mu \in P_{2}(\mathbb{R}^{d})$, we let $\mu_{N}$ be the canonical $N$-particle approximation in Definition~\ref{def-7}. Under the condition~\eqref{2-1-3-1}, we have $(\mu_N)_N$ converges to $\mu$ in $P_2(\mathbb{R}^d)$
\end{lemma}

\begin{proof}
Let
$ X(\theta):=\nabla \phi_{\to\mu}(\theta)$ with $\theta\in[0,1]^d$ be the optimal transport map, where $\phi_{\to\mu}$ is the Brenier potential transporting $\nu_{[0,1]^d}$ to $\mu$. Hence $X_{\#}\nu_{[0,1]^d}=\mu$. Moreover, because $\mu\in P_2(\mathbb{R}^d)$, we also have
\begin{equation*}
    \int_{[0,1]^d}|X(\theta)|^2 \ud \nu_{[0,1]^d}(\theta)
    =\int_{\mathbb{R}^d}|x|^2 \ud \mu(x)<\infty.
\end{equation*}
Thus, $X\in L^2\bigl([0,1]^d,\nu_{[0,1]^d};\mathbb{R}^d\bigr)$.

Consider the piecewise constant map $X_N:[0,1]^d\to\mathbb{R}^d$ defined by
\begin{equation*}
    X_N(\theta)=x_\ell,
    \qquad \theta\in E_\ell,
\end{equation*}
where $x_\ell$ is defined in \eqref{eqn:xl}.
Since $\nu_{[0,1]^d}(E_\ell)=1/N$, it follows that
\begin{equation*}
    (X_N)_{\#}\nu_{[0,1]^d}
    =\frac1N\sum_{\ell=1}^N\delta[x_\ell]
    =\mu_N.
\end{equation*}
Therefore, $(X_N,X)_{\#}\nu_{[0,1]^d}$
is a coupling between $\mu_N$ and $\mu$. Hence by \citep[Equation (7.1.6)]{ambrosio2008gradient},
\begin{equation}\label{eq:canonical-particle-W2-bound}
    W_2(\mu_N,\mu)^2
    \le
    \int_{[0,1]^d}|X_N(\theta)-X(\theta)|^2 \ud \nu_{[0,1]^d}(\theta).
\end{equation}
We next show that $X_N\to X$ in $L^2$. To this end, we let $P_N$ denote the cell-average operator associated with the partition $(E_\ell)_{\ell=1}^N$, defined by
\begin{equation*}
    (P_N\sigma)(\theta)
    :=
    \frac{1}{\nu_{[0,1]^d}(E_\ell)}
    \int_{E_\ell}\sigma(\underline{\theta})\,
    \ud\nu_{[0,1]^d}(\underline{\theta}),
    \qquad \theta\in E_\ell.
\end{equation*}
By construction, one has $X_N=P_NX$. We observe that $P_N$ is an
$L^2$-contraction: $ \|P_N\sigma\|_{L^2}\leq \|\sigma\|_{L^2}$.
Next, for any $Y\in C([0,1]^d;\mathbb R^d)$, the definition of $P_N$
gives
\[
    \|P_NY-Y\|_{L^\infty}
    \leq \omega_Y(h_N)\to0
    \qquad\text{as }~~ h_N\to0,
\]
where $h_N = \max_{1\leq l \leq N} \text{diam} (E_l)$, and $\omega_Y$ denotes the modulus of continuity of $Y$. Since
$C([0,1]^d;\mathbb R^d)$ is dense in
$L^2([0,1]^d;\mathbb R^d)$ and $P_N$ is uniformly bounded on $L^2$,
the convergence extends to every $X\in L^2([0,1]^d;\mathbb R^d)$.
Therefore,
\[
    \|X_N-X\|_{L^2}
    =
    \|P_NX-X\|_{L^2}
    \to0.
\]
Combining this with \eqref{eq:canonical-particle-W2-bound}, we obtain the converge of $(\mu_N)_N$ to $\mu$ in $P_2(\mathbb{R}^d)$, as $N$ tends to $+\infty.$
\end{proof}

\subsection{Functionals in Wasserstein spaces}
In this part, we collect several elementary properties on functionals~\eqref{FVW} and~\eqref{FU} in Wasserstein space $P_{2}(\mathbb{R}^{d})$.
For the internal energy~\eqref{FU}, we let $\rho = \mathrm{d}\mu/\mathrm{d}x$ be the Radon–Nikodym derivative of $\mu$ with respect to the Lebesgue measure, when $m \ge 1$, we define
\begin{equation*}
    \mathcal{F}_{U}(\mu) :=
    \begin{cases}
        \int_{\mathbb{R}^{d}} U(\rho(x)) \ud x &\text{ if } \mu \in P_{2}^{\mathrm{ac}}(\mathbb{R}^{d})\\
        +\infty &\text{ otherwise.}
    \end{cases}
\end{equation*}
and when $m < 1$, we define
\begin{equation*}
    \mathcal{F}_{U}(\mu) :=
        \int_{\mathbb{R}^{d}} U(\rho(x)) \ud x.
\end{equation*}
This difference for $m \ge 1$ and $m < 1$ depends on the recession function $U_{\infty} := \lim_{t \to \infty} U(t)/t$. In particular, for $m \ge 1$, we have $U_{\infty} = +\infty$, and for $m < 1$, we have $U_{\infty} = 0$.

We recall the following lower-semi-continuity of functionals considered here.
\begin{lemma}\label{lemma-7-stability}For any $\mu \in P_{2}(\mathbb{R}^{d})$ and any sequence of $\{\mu_{N}\}_{N = 1}^{\infty}$ in $P_{2}(\mathbb{R}^{d})$ with uniformly bounded second moment, such that $\mu_{N} \rightharpoonup \mu$ in the narrow topology, we have
\begin{equation*}
    \liminf_{N \to \infty} \mathcal{F}_{U}(\mu_{N}) \ge \mathcal{F}_{U}(\mu),\quad \liminf_{N \to \infty} \mathcal{F}_{V}(\mu_{N}) \ge \mathcal{F}_{V}(\mu),\quad \liminf_{N \to \infty} \mathcal{F}_{W}(\mu_{N}) \ge \mathcal{F}_{W}(\mu).
\end{equation*}
and $\liminf_{N \to \infty}W_{2}(\mu_{N}, \nu) \ge W_{2}(\mu, \nu)$.
\end{lemma}
\begin{proof}
The lower-semi-continuity of $\mathcal{F}_{U}$ is given in \cite[Proposition 2.6]{coudreuse2026aronson}, the lower-semi-continuity of $\mathcal{F}_{V}, \mathcal{F}_{W}$ is due to \cite[Lemma 4.3]{villani2008optimal}, and the lower-semi-continuity of $W_{2}$ is given in \cite[Remark 6.12]{villani2008optimal}.
\end{proof}

\begin{lemma}[Proposition 2.7 of \cite{coudreuse2026aronson}]\label{lemma-5-internal}
For $m>m_0:=d/(d+2)$, the internal energy $\mathcal{F}_U$ satisfies the lower bound
\begin{equation*}
    \mathcal{F}_{U}(\mu)
    \geq
    \begin{cases}
        C(1,d)-\dfrac{d}{2}\log M_2(\mu),
        & m=1,\\[2mm]
        C(m,d)M_2(\mu)^{\frac{d}{2}(1-m)},
        & m\neq1.
    \end{cases}
\end{equation*}
In particular, $\mathcal{F}_{U}$ is well defined as an extended-real-valued functional $\mathcal{F}_{U}:P_2(\mathbb{R}^d)\to(-\infty,+\infty]$,
in the sense that $\mathcal{F}_{U}(\mu)>-\infty$ for every $\mu\in P_2(\mathbb{R}^d)$. Moreover, since $\frac{d}{2}(1-m)<1$ for all $m>m_0$, 
its negative part grows at most sublinearly with respect to the second moment.
\end{lemma}

\begin{lemma}\label{lemma-6-well} 
Suppose that, for some constants $C>0$ and $0<\beta<2$,
\begin{equation*}
    V(x)\geq -C(1+|x|^\beta),
    \qquad
    W(x)\geq -C(1+|x|^\beta).
\end{equation*}
Then, for every $\mu\in P_2(\mathbb{R}^d)$,
\begin{equation*}
    \mathcal{F}_{V}(\mu)
    \geq
    -C\big(1+M_2(\mu)^{\beta/2}\big), \quad 
        \mathcal{F}_{W}(\mu)
    \geq
    -2C\big(1+M_2(\mu)^{\beta/2}\big)\,.
\end{equation*}
Consequently, both $ \mathcal{F}_{V},\mathcal{F}_{W}:
    P_2(\mathbb{R}^d)\to(-\infty,+\infty]$
are well defined as extended-real-valued functionals. Moreover, they are uniformly bounded from below on every subset of $P_2(\mathbb{R}^d)$ with uniformly bounded second moments.
\end{lemma}

\begin{proof}
The proof is a direct result of H\"older's inequality.
\end{proof}

\begin{lemma}[coercivity]\label{lemma-11-coercive}Given $\nu \in P_{2}(\mathbb{R}^{d})$, the JKO functional is coercive in $P_{2}(\mathbb{R}^{d})$ in the sense that
\begin{equation*}
    \JKO(\mu,\nu) \ge \frac{1}{4\tau}M_{2}(\mu) - \frac{1}{2\tau}M_{2}(\nu) - C\big(1 + M_{2}(\mu)^{\beta/2}\big)
\end{equation*}
for some positive constant $\beta < 2$,  and $C$ is a constant depending on $U, V, W$. Hence, for every $E\in\mathbb{R}$, the sublevel set $\{\mu: \JKO(\mu, \nu) \le E\}$ has uniformly bounded second moments.
\end{lemma}
\begin{proof}
Combining Lemmas~\ref{lemma-5-internal}, and~\ref{lemma-6-well}, we get
\begin{equation*}
    \mathcal{F}(\mu) = \mathcal{F}_U(\mu) + \mathcal{F}_V(\mu) + \mathcal{F}_W(\mu) \ge -C(1 + M_{2}(\mu)^{\beta/2})
\end{equation*}
for positive constants $C$ and $\beta < 2$. On the other hand, we apply the triangle inequality in \cite[Equations (17) and (18)]{JKO} to get 
\begin{equation*}
    W_{2}(\mu, \nu)^{2} \ge \frac{1}{2}M_{2}(\mu) - M_{2}(\nu),
\end{equation*}
which proves the inequality.
\end{proof}
We next present an important observation by Coudreuse \citep[Section 3.3]{coudreuse2026aronson}.
\begin{lemma} \label{lemma-12}
Let $\nu,\nu_N\in P_2(\mathbb R^d)$ satisfy $\nu_N \to \nu~\text{in } P_2(\mathbb R^d).$ Suppose that that
\begin{equation*}
    \Gamma\text{-}\lim_{N \to \infty}\JKO_N(\cdot,\nu_N) = \JKO(\cdot,\nu), \quad \Gamma\text{-}\lim_{N \to \infty} \mathcal F_N = \mathcal F,
\end{equation*}
where the $\Gamma$-convergence is understood along sequences $\mu_N \rightharpoonup \mu$ narrowly satisfying a uniform second moment bound $\sup_{N} M_2(\mu_N)<\infty$. Assume there are constants $C>0$ and $0<\beta<2$ such that
\begin{equation}\label{eq:nov11.2026.3}
\mathcal F_N(\mu_N)\ge-C\bigl(1+M_2(\mu_N)^{\beta/2}\bigr)
\end{equation}
for every $N$ and every admissible $\mu_N$. Let $\mu_N$ be a minimizer of $\JKO_N(\cdot,\nu_N)$, and $\mu_*$ be the unique minimizer of $\JKO(\cdot,\nu)$. Then $(\mu_N)_N$ converges to $\mu_*$ in $P_2(\mathbb R^d)$.
\end{lemma}
\begin{proof}
Set
\[
\mathcal J_N(\mu):=\JKO_N(\mu,\nu_N)=\frac{1}{2\tau}W_2(\mu,\nu_N)^2+\mathcal F_N(\mu), \quad \mathcal J(\mu):=\JKO(\mu,\nu)=\frac{1}{2\tau} W_2(\mu,\nu)^2 +\mathcal F(\mu).
\]
Using \eqref{eq:nov11.2026.3}, 
\[
\mathcal J_N(\mu_N) =\mathcal F_N(\mu_N) +\frac1{2\tau}W_2^2(\mu_N,\nu_N) \ge \frac1{4\tau}M_2(\mu_N)-\frac1{2\tau}M_2(\nu_N)-C\bigl(1+M_2(\mu_N)^{\beta/2}\bigr).
\]
Because $\nu_N\to\nu$ in $W_2$, $\big(M_2(\nu_N)\big)_N$ is a bounded sequence in $\mathbb R.$ Since $\beta<2$, the sublevel sets of $\mathcal J_N$ have uniformly bounded second moments. Let $(\bar \mu_N)$ be a recovery sequence for $\mu^*$. By minimality,
\[
\mathcal J_N(\mu_N)\le \mathcal J_N(\bar \mu_N)\rightarrow \mathcal J(\mu^*) \]
Hence $(\mu_N)$ has uniformly bounded second moments and is relatively compact in the narrow topology.

Building on the $\Gamma$-convergence of the JKO functional, and the uniqueness of the minimizer of the JKO functional $\mathcal{J}$, we can apply the fundamental theorem of $\Gamma$-convergence to get the following (see, for example, \cite[Theorem 1.21]{braides2002gamma} and \cite[Theorem 2.1]{braides2014local}):
\[
    \mu_N \rightharpoonup \mu_* ~~ \text{narrowly}, \quad \text{and} \quad \mathcal J_N(\mu_N) \rightarrow \mathcal J(\mu_\ast). 
\]
Now if we let
\[
    \alpha_0 := \frac{1}{2\tau} W_2(\mu_*,\nu)^2,
    \quad
    \alpha_1 := \mathcal{F}(\mu_*).
\]
then we can split the convergence of energy as
\[
    \mathcal{J}_N(\mu_N) \to \alpha_0+\alpha_1.
\]
By the lower semicontinuity of the Wasserstein distance and the $\Gamma$-liminf inequality for $\mathcal{F}_N$, we have
\begin{equation}\label{2-9-wasserstein}
    \liminf_{N\to\infty}
    \frac{1}{2\tau}W_2(\mu_N,\nu_N)^2
    \geq \alpha_0,
    \qquad
    \liminf_{N\to\infty}
    \mathcal{F}_N(\mu_N)
    \geq \alpha_1.
\end{equation}
Consequently,
\begin{align*}
\limsup_{N\to\infty}  \frac{1}{2\tau}W_2(\mu_N,\nu_N)^2\leq \lim_{N\to\infty}\mathcal{J}_N(\mu_N)  - \liminf_{N\to\infty}\mathcal{F}_N(\mu_N) \leq (\alpha_0+\alpha_1)-\alpha_1  =\alpha_0.
\end{align*}
Combining the lower semi-continuity of the Wasserstein distance in~\eqref{2-9-wasserstein},
we can get the convergence of the distance:
\[
\lim_{N \to \infty} W_2(\mu_N,\nu_N) = W_2(\mu_*,\nu).
\]

Building on the convergence of distance, Coudreuse \cite[Section 3.3]{coudreuse2025li} considered the following identity of moment to upgrade the narrow convergence $\mu_{N} \rightharpoonup \mu$ to a $P_{2}(\mathbb{R}^{d})$ convergence,
\begin{equation}\label{2-8-extra}
    M_2(\mu_N)=2M_2(\nu_N)+2W_2(\mu_N,\nu_N)^2-\int_{\mathbb R^d\times\mathbb R^d}|x-2y|^2 \ud \gamma_{\mathrm{o}, N}(x,y)
\end{equation}
where $\gamma_{\mathrm{o}, N}$ is an optimal coupling between $\mu_N$ and $\nu_N$ under the 2-Wasserstein distance. In particular, since $\mu_N\to\mu_*$,
$\nu_N\to\nu$ narrowly, one has $\gamma_{\mathrm{o}, N}\to\gamma_{\mathrm{o}}$ narrow according to \cite[Proposition 7.1.3]{ambrosio2008gradient}, where $\gamma_{\mathrm{o}}$ is an optimal coupling between $\mu_*$ and $\nu$ under the 2-Wasserstein distance.

Using \cite[Lemma 4.3]{villani2008optimal}, we have
\[
\liminf_{N\to\infty}
\int_{\mathbb R^d\times\mathbb R^d} |x-2y|^2\ud\gamma_{\mathrm{o}, N}
\geq
\int_{\mathbb R^d\times\mathbb R^d} |x-2y|^2\ud\gamma_{\mathrm{o}}.
\]
Since $M_2(\nu_N)\to M_2(\nu)$ and $W_2(\mu_N,\nu_N)\to W_2(\mu_*,\nu)$, 
we obtain
\[
\begin{aligned}
\limsup_{N\to\infty}M_2(\mu_N)
\leq
2M_2(\nu)
+2W_2(\mu_*,\nu)^2
-\int_{\mathbb R^d\times\mathbb R^d} |x-2y|^2\ud\gamma_{\mathrm{o}} = 
M_2(\mu_*).
\end{aligned}
\]
On the other hand, the second moment is lower-semicontinous under narrow convergence:
\begin{equation*}
    \liminf_{N\to\infty}M_2(\mu_N) \ge M_2(\mu_*)\,.
\end{equation*} 
Hence, $M_2(\mu_N)\to M_2(\mu_*)$. Since $\mu_N\to\mu_*$ narrowly, we infer $\mu_N \to \mu_*$ in $P_{2}(\mathbb{R}^{d})$.  
\end{proof}
\begin{remark}This upgrade allows the convergence of minimizers to be propagated through successive JKO time steps: the convergence in $P_2(\mathbb R^d)$ at one time step provides the hypothesis needed to establish convergence at the next. This is nontrivial because $P_2(\mathbb R^d)$ is not precompact in the $W_2$ metric topology. The JKO sub-level bounds in Lemma~\ref{lemma-11-coercive} restrict the minimizing sequence to a family with uniformly bounded second moments, which yields the precompactness of discrete JKO minimizers in the narrow topology and $P_{p}(\mathbb{R}^{d})$ topology where $p < 2$, but not in $W_2$ metric topology (see, for example, \cite[Remark 2.22]{di2024optimal}). The key additional ingredient is the convergence $W_2(\mu_N,\nu_N) \to W_2(\mu_*,\nu)$, which, together with $\nu_N\to\nu$ in $P_2(\mathbb R^d)$ yields convergence of the second moments via the identity~\eqref{2-8-extra}. This upgrades the narrow convergence of $\mu_N$ to convergence in the $W_2$ metric.
\end{remark}

\section{Convergence on real line}\label{sec:3}
In this section, we establish the convergence of the particle methods for JKO schemes on $\mathbb{R}$ shown in Theorem~\ref{theorem:real}. The particle method is defined by minimizing the following discrete JKO functional.
\begin{definition}[discrete JKO functional, one dimension]\rm Given $N \in \mathbb{N}$ and $\nu_N = \frac{1}{N} \sum_{\ell=1}^N \delta[y_\ell]$ with $y_1 \leq  y_2 \leq \cdots \leq y_N $, we define
\begin{align} \label{pJKO-1D}
    \JKO_{N}(\mu_{N}; \nu_N) &:= \mathcal{F}_{N}(\mu_{N}) + \frac{1}{2\tau}W_{2}\big(\mu_{N}, \nu\big)^{2}
    = \mathcal{F}_{N}(\mu_{N}) + \frac{1}{2 N\tau} \sum_{\ell = 1}^N |x_\ell - y_\ell|^2
\end{align}
where $x_1 \leq x_2 \cdots \leq x_N$  is defined in \eqref{eqn:xl} with $d=1$, and $\mu_{N} = \frac{1}{N}\sum_{\ell = 1}^{N} \delta[x_{\ell}]$ is the empirical measure associated with the particles. $X_{\nu}$ is the inverse distribution of $\nu$. Here, the discrete free energy is defined by
\begin{equation*}
    \mathcal{F}_{N}(\mu_{N}) := \mathcal{F}_{U, N}(\mu_{N}) + \mathcal{F}_{V}(\mu_{N}) + \mathcal{F}_{W}(\mu_{N})\,,
\end{equation*}
where
\[
 \mathcal{F}_{V}(\mu_{N}) = \frac{1}{N} \sum_{\ell=1}^N V(x_\ell), \quad \mathcal{F}_{W}(\mu_{N}) = \frac{1}{2N^2} \sum_{j, \ell = 1}^N W(x_j - x_\ell)\,,
\]
and 
\begin{equation}\label{3-9-internal}
    \mathcal{F}_{U, N}(\mu_{N}) := \sum_{\ell = 1}^{N - 1} \frac{1}{N} U_{\mathrm{mc}} \big(N(x_{\ell + 1} - x_{\ell})\big)\,.
\end{equation}
Here $U_{\mathrm{mc}}$ is a lower-semi-continuous function given by
\begin{equation*}
    U_{\mathrm{mc}}(\alpha) = \begin{cases}
        \alpha U\big(\frac{1}{\alpha}\big) &\text{ when } \alpha > 0\\
        U_{\infty} &\text{ when } \alpha = 0\\
        +\infty &\text{ when } \alpha < 0\,,
    \end{cases}
\end{equation*}
and $U_{\infty} := \lim_{\alpha \rightarrow 0} U_\text{mc}(\alpha)=  \lim_{t \to \infty} \frac{U(t)}{t}$ is the \textit{recession function} of $U$.
\end{definition}

\subsection{Existence of minimizer}
\begin{theorem}[existence]\label{3-well-pose}
Given $N \in \mathbb{N}$, and $\nu_N \in P_{2, N}(\mathbb{R})$, there exists an $N$-particle solution $\mu_{N}$ to the discrete JKO problem \eqref{pJKO-1D}. 
\end{theorem}

To prove this theorem, we define an auxiliary measure $\widetilde\mu_N$ to replace $\mu_{N} = \frac{1}{N}\sum_{\ell = 1}^{N} \delta[x_{\ell}]$ for the  internal energy. Set $\theta_\ell = (\ell-\frac12)/N$ and 
\begin{equation}\label{eq:probability-reconstruction}
X_N(\theta)=
\begin{cases}
  x_1+\theta-\theta_1,
     & 0\leq\theta\leq\theta_1,\\[2pt]
  x_\ell + N (\theta-\theta_\ell) (x_{\ell+1}-x_\ell),
     & \theta_\ell\leq\theta\leq\theta_{\ell+1},\quad 1\leq \ell<N,\\[2pt]
  x_N+\theta-\theta_N,
     & \theta_N\leq\theta\leq1.
\end{cases}
\end{equation}
If $N=1$, we take $X_1(\theta)=x_1+\theta-\frac12$.
We then define the following reconstructed probability measure (or an auxiliary measure)
\[
  \widetilde\mu_N:=(X_N)_\#\bigl(\mathcal L^1|_{[0,1]}\bigr)
    \in P_2(\mathbb R).
\]
Here $\mathcal L^1|_{[0,1]}$ is the one-dimensional Lebesgue measure restricted to $[0, 1]$ which is also the uniform distribution. We note that, if $x_1< x_2< \cdots <x_N$, the auxiliary measure $\widetilde\mu_N$ does not possess atoms and  
\begin{equation}\label{eq:reconstruction-energy}
 \mathcal F_U(\widetilde\mu_N)
   =\int_0^1 U_{\mathrm{mc}}(X_N'(\theta)) \ud \theta
   =\mathcal F_{U,N}(\mu_N)+\frac{U_{\mathrm{mc}}(1)}N = \mathcal F_{U,N}(\mu_N)+\frac{U(1)}N.
\end{equation}
One checks that 
\begin{equation}\label{eq:reconstruction-moment-bounds} 
 \frac23 M_2(\mu_N)
 \leq M_2(\widetilde\mu_N)
 \leq M_2(\mu_N)
 +\frac{\sqrt{2M_2(\mu_N)}}{4N^{3/2}}
 +\frac1{12N^3}.
\end{equation}
Also, for the cumulative distribution functions, 
\begin{equation}\label{eq:reconstruction-cdf}
 \sup_{t\in\mathbb R}\left|\widetilde\mu_N((\! -\infty,t])              -\mu_N((\! -\infty,t])\right| \leq \frac1{2N}.
\end{equation}

Using this reconstructed probability measure, we get the following two lemmas.
\begin{lemma} \label{lem9}
There exist $C>0$ and $0<\beta<2$, independent of $N$, such that for any $\mu_{N} \in P_{2, N}(\mathbb{R})$, we have,
\begin{equation*}
    \mathcal{F}_{U, N}(\mu_{N}) \ge -C\big(1 + M_{2}(\mu_{N})^{\beta/2}\big).
\end{equation*}
\end{lemma}

\begin{proof} By using~\eqref{eq:reconstruction-energy}, we have
\[
\mathcal F_{U,N}(\mu_N)=\mathcal F_U(\widetilde\mu_N)-\frac{U(1)}N.
\]
Then, applying Lemma~\ref{lemma-5-internal}, there exist $C>0$ and $0<\beta<2$, independent of $N$, such that
\[
\mathcal F_U(\widetilde\mu_N) \ge -C\big(1+M_2(\widetilde\mu_N)^{\beta/2}\big).
\]
By using \eqref{eq:reconstruction-moment-bounds}, there exists a constant $C_1$ independent of $N$ such that
\[
M_2(\widetilde\mu_N)\le M_2(\mu_N)+\frac{\sqrt{2M_2(\mu_N)}}{4N^{3/2}}+\frac1{12N^3}\le C_1(1+M_2(\mu_N)),
\]
and consequently, we have
\[
\mathcal F_{U,N}(\mu_N)\ge -C_2\big(1+M_2(\mu_N)^{\beta/2}\big),
\]
for a constant $C_2>0$ independent of $N$. This proves the Lemma.
\end{proof}

\begin{lemma}\label{lemma5} Let $\bx :=(x_1, \cdots, x_N)^T$. For all $(\bx_{(j)})_{j = 1}^{\infty} \to \bx$ in $\mathbb{R}^{N}$, we have
\begin{equation*}
    \liminf_{j \to \infty} \mathcal{F}_{U, N}(\mu_{j}) \ge \mathcal{F}_{U, N}(\mu)\,,
\end{equation*}
where $\mu_{j} = \frac{1}{N}\sum_{\ell = 1}^{N} \delta[x_{(j), \ell}]$ and $\mu = \frac{1}{N} \sum_{\ell = 1}^{N} \delta[x_{\ell}]$.
\end{lemma}

\begin{proof} Let $X_N^{(j)}$ and $X_N$ be the piecewise-affine maps defined by \eqref{eq:probability-reconstruction}, corresponding respectively to $x^{(j)}$ and $x$. Set
\[
\widetilde\mu_N^{(j)}=(X_N^{(j)})_\#(\mathcal L^1|_{[0,1]}),\quad\widetilde\mu_N=(X_N)_\#(\mathcal L^1|_{[0,1]}).
\]
Since $x^{(j)}\to x$ in $\mathbb R^N$, we have $X_N^{(j)}\to X_N$ uniformly on $[0,1]$. Therefore,
\[
W_2(\widetilde\mu_N^{(j)},\widetilde\mu_N)^2 \le \int_0^1|X_N^{(j)}(\theta)-X_N(\theta)|^2 \ud \theta \rightarrow0 \text{ as } N \to \infty.
\]
By using Lemma \ref{lemma-7-stability} and identity \eqref{eq:reconstruction-energy},
\[
\liminf_{j\to\infty}\mathcal F_{U,N}(\mu^{(j)})=\liminf_{j\to\infty}\left(\mathcal F_U(\widetilde\mu_N^{(j)})-\frac{U(1)}N\right)
\ge \mathcal F_U(\widetilde\mu_N)-\frac{U(1)}N=\mathcal F_{U,N}(\mu).
\]
This proves the Lemma.
\end{proof}

\begin{proof}[Proof of Theorem~\ref{3-well-pose}]
According to Lemma~\ref{lem9} and Lemma~\ref{lemma-6-well}, we have
\begin{equation*}
    \mathcal{F}_{N}(\mu_{N}) = \mathcal{F}_{U, N}(\mu_{N}) + \mathcal{F}_{V}(\mu_{N}) + \mathcal{F}_{W}(\mu_{N}) \ge -C(1 + M_{2}(\mu_{N})^{\beta/2})
\end{equation*}
for all $\mu_{N} \in P_{2, N}(\mathbb{R}^{d})$ with $\beta < 2$. Hence, by triangle inequality, we get
\begin{equation}\label{3-1-JKO}
    \JKO_{N}(\mu_{N}, \nu_N) = \mathcal{F}_{N}(\mu_{N}) + \frac{1}{2\tau} W_{2}\big(\mu_{N}, \nu_N\big)^{2} \ge \frac{1}{4\tau} M_{2}(\mu_{N}) - \frac{1}{2\tau} M_{2}(\nu_N) -C(1 + M_{2}(\mu_{N})^{\beta/2})
\end{equation}
which means that $\JKO_{N}(\cdot, \nu_N)$ is bounded from below on $P_{2, N}(\mathbb{R})$.

Now, let $(\mu_{j})_{j} = (\frac{1}{N}\sum_{\ell = 1}^{N} \delta[x_{(j), \ell}])_{j}$ be a minimizing sequence for $\JKO_{N}(\cdot, \nu_N)$. By the coercivity estimate in \eqref{3-1-JKO}, there exists a constant $C>0$, independent of $j$, such that
\begin{equation*}
    \sup_j \|\bx_{(j)}\|^{2} = N \times \sup_j M_{2}(\mu_{j}) \le C\,.
\end{equation*}
Since $(\bx_{(j)})_j$ is bounded in the finite-dimensional space $\mathbb{R}^N$, it admits a convergent subsequence by Bolzano-Weierstrass theorem. Thus, up to extracting a subsequence, there exists $\bx \in \mathbb{R}^N$ such that $(\bx_{(j)})_j$ converges to $\bx$ in $\mathbb{R}^N$. Furthermore, the set
\begin{equation*}
    \left\{
    \bx \in \mathbb{R}^N :
    x_\ell \le x_{\ell+1},
    \ \ell = 1,\ldots,N-1
    \right\}
\end{equation*}
is closed. Hence the ordering is preserved in the limit, and we have
\begin{equation*}
    x_\ell \le x_{\ell+1},
    \quad \ell = 1,\ldots,N-1.
\end{equation*}

We associate this limiting particle position with an empirical particle measure defined as follows,
\begin{equation*}
    \mu_{\ast} := \frac{1}{N} \sum_{\ell = 1}^{N} \delta[x_{\ell}] \in P_{2, N}(\mathbb{R}^{d}).
\end{equation*}
We now claim that $\mu_*$ is the desired minimizer. Indeed, by triangle inequality, we have
\begin{equation*}
    W_{2}(\mu_{j}, \nu_N) - W_{2}(\mu_{\ast}, \nu_N) \le W_{2}(\mu_{j}, \mu_{\ast}) = \left(\frac1N\sum_{\ell=1}^N|x_{(j),\ell}-x_\ell|^2\right)^{1/2}.
\end{equation*}
Hence, $\lim_{j \to \infty} W_{2}(\mu_{j}, \nu_N) = W_{2}(\mu_{\ast}, \nu_N)$. By the lower semi-continuity of the free energies (Lemma~\ref{lemma5} and Lemma~\ref{lemma-7-stability}), we obtain
\begin{equation*}
    \mathcal{F}_{N}(\mu_{\ast}) \le \liminf_{j \to \infty} \mathcal{F}_{U, N}(\mu_{j}) + \liminf_{j \to \infty} \, (\mathcal{F}_{V}(\mu_{j}) + \mathcal{F}_{W}(\mu_{j})) =: \liminf_{j \to \infty} \mathcal{F}_{N}(\mu_{j}).
\end{equation*}
Hence,
\begin{equation*}
    \JKO_{N}(\mu_{\ast}, \nu_N) \le \liminf_{j \to \infty} \JKO_{N}(\mu_{j}, \nu_N) = \min_{P_{2}(\mathbb{R})} \JKO_{N}(\cdot, \nu_N)\,,
\end{equation*}
which means that $\mu_{*}$ is the desired minimal.
\end{proof}

\subsection{Convergence in large sample limit}
In this part, we establish the convergence of the particle method in Theorem~\ref{theorem:real}. In view of Lemma~\ref{lemma-12}, it suffices to prove the $\Gamma$-convergence of both the energy functional and the JKO loss. The convergence of the particle method then follows by a mathematical induction argument.

\begin{theorem}[$\Gamma$-convergence]\label{theorem6-gamma} Given $(\nu_{N})_{N = 1}^{\infty}$  converges to $\nu$ in $P_{2}(\mathbb{R})$, we have 
\begin{equation*}
     \Gamma\text{-}\lim_{N \to \infty} \JKO_{N}(\cdot, \nu_{N}) = \JKO(\cdot, \nu), \quad  \Gamma\text{-}\lim_{N \to \infty} \mathcal F_N  = \mathcal F.
\end{equation*}
The following assertions hold, for any $\mu \in P_{2}(\mathbb{R})$, and any $\mu_{N} \in P_{2, N}(\mathbb{R})$ such that $\mu_{N} \to \mu$ narrowly with uniformly bounded second moment:
\begin{equation*}
    \liminf_{N \to \infty} \JKO_{N}(\mu_{N}, \nu_{N}) \ge \JKO(\mu, \nu), \quad  
    \liminf_{N \to \infty} \mathcal F_{N}(\mu_{N}) \ge \mathcal F(\mu)\,.
\end{equation*}
On the other hand, for any $\mu \in P_{2}(\mathbb{R})$, we let $\mu_{N}$ be its $N$-canonical sampling in Definition~\ref{def-7}, we have
\begin{equation*}
    \limsup_{N \to \infty} \JKO_{N}(\mu_{N}, \nu_{N}) \le \JKO(\mu, \nu), \quad 
       \limsup_{N \to \infty} \mathcal F_{N}(\mu_{N}) \le \mathcal F(\mu)\,.
\end{equation*}
\end{theorem}
\begin{proof} ($\Gamma$-$\liminf$ inequality) Let $\mu_N \in  P_{2,N}(\mathbb R)$ be such that $\big(\mu_N\big)_N$ converging narrowly to $\mu$ and $(M_2(\mu_N))_N$ is uniformly bounded in $\mathbb R$. Let $\widetilde\mu_N$ be the reconstructed probability measure defined by \eqref{eq:probability-reconstruction}. By~\eqref{eq:reconstruction-cdf}, $(\widetilde\mu_N)_N$ converges narrowly to $\mu$. Furthermore, \eqref{eq:reconstruction-moment-bounds}  implies $(M_2(\widetilde\mu_N))_N$ is bounded in $\mathbb R.$ Applying Lemma \ref{lemma-7-stability} and using \eqref{eq:reconstruction-energy}, we obtain
\[
\liminf_{N\to\infty}\mathcal F_{U,N}(\mu_N)=\liminf_{N\to\infty}\left(\mathcal F_U(\widetilde\mu_N)-\frac{U(1)}N\right)\ge \mathcal F_U(\mu).
\]
Again by  Lemma \ref{lemma-7-stability}, the confinement and interaction energies satisfy
\[
\liminf_{N\to\infty}\mathcal F_V(\mu_N) \ge\mathcal F_V(\mu), \qquad \liminf_{N\to\infty}\mathcal F_W(\mu_N)\ge\mathcal F_W(\mu).
\]
Therefore,
\[
\liminf_{N\to\infty}\mathcal F_N(\mu_N)\ge\mathcal F(\mu).
\]
Finally, since $\nu_N\to\nu$ in $W_2$, 
\[
\liminf_{N\to\infty}W_2^2(\mu_N,\nu_N) \ge W_2^2(\mu,\nu), \quad \text{and hence} \quad \liminf_{N\to\infty}\JKO_N(\mu_N,\nu_N)\ge \JKO(\mu,\nu).
\]
This proves the $\Gamma$-liminf inequality.

~\\
We claim that the $N$-canonical sampling defined in Definition~\ref{def-7} is the recovery sequence for each $\mu \in P_{2}(\mathbb{R})$. 
Indeed, assume $\mathcal F_U(\mu)<\infty$; then $\mathcal F_U(\mu)=\int_0^1 U_{\mathrm{mc}}(X_\mu'(s)) \ud s$, with the same recession convention as above.
Set $h=1/N$ and $x_\ell=h^{-1}\int_{(\ell-1)h}^{\ell h}X_\mu(\theta) \ud \theta$.  For $I_\ell=((\ell-1)h, \ell h)$,
\[ N(x_{\ell+1}-x_\ell)=N^2\int_{I_\ell}[X_\mu(\theta+h)-X_\mu(\theta)] \ud \theta\geq N^2\int_{I_\ell}\int_0^h X_\mu'(\theta+s) \ud s \ud \theta.\]
Since $m>d/(d+2)>1-1/d$ when $d=1$, $U_{\mathrm{mc}}$ is convex and nonincreasing on $(0,\infty)$ and so, Jensen's inequality gives

\begin{equation*}
    \mathcal{F}_{U, N}(\mu_{N}) = \sum_{\ell = 1}^{N - 1}\frac{1}{N} U_{\mathrm{mc}}\big(N(x_{\ell + 1} - x_{\ell})\big) \le \sum_{\ell = 1}^{N - 1}\int_{(\ell - 1)/N}^{\ell/N} \bigg[N \int_{0}^{1/N} U_{\mathrm{mc}}\big(X_{\mu}^{\prime}(\theta + s)\big) \ud s\bigg] \ud \theta.
\end{equation*}
The right hand side has a limit
\begin{equation*}
    \sum_{\ell = 1}^{N - 1}\int_{(\ell - 1)/N}^{\ell/N} \bigg[N \int_{0}^{1/N} U_{\mathrm{mc}}\big(X_{\mu}^{\prime}(\theta + s)\big) \ud s\bigg] \ud \theta \to \int_{0}^{1}U_{\mathrm{mc}}\big(X_{\mu}^{\prime}(\theta)\big) \ud \theta = \mathcal{F}_{U}(\mu) \text{ as $N \to \infty$.} 
\end{equation*}
Hence, we have the $\limsup$ inequality:
\begin{equation*}
    \limsup_{N \to \infty} \mathcal{F}_{U, N}(\mu_{N}) \le \mathcal{F}_{U}(\mu)\,.
\end{equation*}
For the confinement energy $\mathcal{F}_{V}$ and the interaction energy
$\mathcal{F}_{W}$, we refer to Proposition~\ref{lemma-17} below, where it is
shown that
\begin{equation*}
    \limsup_{N\to\infty}\mathcal{F}_{V}(\mu_N)
    \leq \mathcal{F}_{V}(\mu),
    \qquad
    \limsup_{N\to\infty}\mathcal{F}_{W}(\mu_N)
    \leq \mathcal{F}_{W}(\mu).
\end{equation*}
Combining these estimates with the corresponding recovery-sequence
estimate for the internal energy, we obtain
\begin{equation*}
    \limsup_{N\to\infty}\mathcal{F}_{N}(\mu_N)
    \leq \mathcal{F}(\mu),
\end{equation*}
where $\mu_N$ denotes the canonical particle approximation of $\mu$.

For the distance part, we use triangle inequality to show that
\begin{equation*}
    \big|W_{2}(\mu_{N}, \nu_{N}) - W_{2}(\mu, \nu)\big| \le W_{2}(\mu_{N}, \mu) + W_{2}(\nu_{N}, \nu) \to 0 \text{ as } N \to \infty,
\end{equation*}
where $W_{2}(\nu_{N}, \nu) \to 0$ comes from the assumption of the theorem, and $W_{2}(\mu_{N}, \mu) \to 0$ due to Lemma~\ref{lemma-7-4}. Hence,
\begin{equation*}
    \limsup_{N \to \infty} \JKO_{N}(\mu_{N}, \nu_{N}) \le \JKO(\mu, \nu).
\end{equation*}
which shows the desired $\Gamma$-convergence.\end{proof}
\begin{proof}[Proof of Theorem~\ref{theorem:real}] We prove the theorem by induction on $k$. 
Let $\mu_{\tk,N}$ denote the $N$-particle solution of the discrete JKO scheme at time step $\tk$. By the coercivity of the JKO functional on $P_{2}(\mathbb{R})$, the sequence $\big(M_{2}(\mu_{\tk,N})\big)_N$ is bounded in $P_{2}(\mathbb{R})$ for each $k.$  
Hence, by Prokhorov's theorem,  a subsequence of $(\mu_{\tk,N})_N$ converges narrowly to some $\mu_{\tk,\ast} \in P_{2}(\mathbb{R}),$ as $N$ tends to $+\infty.$
We next use the $\Gamma$-convergence result in Theorem~\ref{theorem6-gamma} to identify the limit $\mu_{\tk,\ast}$.

Indeed, since the initially $\big(\mu_{0,N}\big)_N$ converges to $\mu_{0}$ in $P_{2}(\mathbb{R})$,  Theorem~\ref{theorem6-gamma}, together with Lemma~\ref{lemma-12}, yields 
$\big(\mu_{1,N}\big)_N$ converges to $\mu_{1}$ in $P_{2}(\mathbb{R}).$
Here $\mu_{1}$ is the first-step solution of the limiting JKO scheme. Since $\mu_{1}$ is unique, every convergent subsequence has the same limit, and hence the entire sequence $\{\mu_{1,N}\}_{N}$ converges to $\mu_{1}$, i.e., $ \mu_{1,\ast}=\mu_{1}$. Iteratively, we obtain that for every $\tk$, $\big(\mu_{\tk,N}\big)_N$ converges to $\mu_{\tk,\ast}=\mu_{\tk}$ in $P_{2}(\mathbb{R}),$ 
where $\mu_{\tk}$ denotes the solution of the limiting JKO scheme at time step $\tk$. \end{proof}

\section{Blob methods}\label{sec:4}
In this section, we investigate the convergence of blob methods. We recall that, at the fully discrete level, the blob method is constructed by minimizing the discrete JKO functional given as follows.

\begin{definition}[discrete JKO functional, blob method] \label{def10}
\rm Given $N \in \mathbb{N}$ and $\nu_N = \frac{1}{N} \sum_{j=1}^N \delta[y_j]$, set 
\begin{equation} \label{JKO-blob}
    \JKO_{N}(\mu_{N}; \nu_N) :=  \mathcal{F}_{\varepsilon, U}(\mu_{N}) + \mathcal{F}_{V}(\mu_{N}) + \mathcal{F}_{W}(\mu_{N}) + \frac{1}{2\tau}W_{2}(\mu_{N}, \nu_N)^{2},
\end{equation}
where $\mu_{N} = \frac{1}{N}\sum_{\ell = 1}^{N} \delta[x_{\ell}]$, and $\mathcal F_V(\mu_N)$, $\mathcal F_W(\mu_N)$, and $\mathcal F_{\varepsilon, U}(\mu_N)$ are defined in \eqref{F_VN} and \eqref{2-bob}, respectively. 
\end{definition}

\subsection{Properties of regularized internal energy}
\begin{lemma}For the internal energies given in Assumption~\ref{assumption}, we have
\begin{equation*}
    \mathcal{F}_{\varepsilon, U}(\mu) \ge -C\big(1 + M_{2}(\mu)^{\beta/2}\big)
\end{equation*}
for some positive constant $C$ and $\beta < 2$.
\end{lemma}

\begin{proof}
By the lower bound in Lemma~\ref{lemma-5-internal}, we obtain
\begin{equation*}
    \mathcal{F}_{\varepsilon, U}(\mu) = \mathcal{F}_{U}(\psi_{\varepsilon} \ast \mu) \ge -C\big(1 + M_{2}(\psi_{\varepsilon} \ast \mu)^{\beta/2}\big)
\end{equation*}
for some positive constant $C$ and $\beta < 2$. Moreover, by the definition of convolution, we have the following estimate for the second moment; see also \cite[Proof of Proposition 3.8]{carrillo2019blob}:
\begin{equation}\label{4-1-1}
    \begin{aligned}
        M_{2}(\psi_{\varepsilon} \ast \mu) &= \int_{\mathbb{R}^{d}} |x|^{2} (\psi_{\varepsilon} \ast \mu)(x) \ud x\\
        &= \int_{\mathbb{R}^{d}} |x|^{2} \int_{\mathbb{R}^{d}} \psi_{\varepsilon}(x - y) \ud \mu(y) \mathrm{d} x = \int_{\mathbb{R}^{d} \times \mathbb{R}^{d}} |x + y|^{2} \psi_{\varepsilon}(x) \ud \mu(y) \mathrm{d} x\\
        &\le \int_{\mathbb{R}^{d} \times \mathbb{R}^{d}} \big(2|x| + 2|y|^{2}\big) \psi_{\varepsilon}(x) \ud \mu(y) \mathrm{d} x = 2\varepsilon^{2} M_{2}(\psi) + 2M_{2}(\mu).
    \end{aligned}
\end{equation}
Absorbing the constant $2\varepsilon^{2} M_{2}(\psi)$ and the remaining scaling factors into the constant $C$, we obtain the desired estimate and complete the proof.
\end{proof}

\begin{lemma}\label{lemma-13-stability}
Fix $\varepsilon>0$. Let $\mu\in P_2(\mathbb{R}^d)$ and let
$(\mu_N)_{N=1}^{\infty}$ be a bounded sequence in $P_2(\mathbb{R}^d),$ and which converges narrowly to $\mu.$
Then
\begin{equation*}
    \liminf_{N\to\infty}
    \mathcal{F}_{\varepsilon,U}(\mu_N)
    \geq
    \mathcal{F}_{\varepsilon,U}(\mu).
\end{equation*}
\end{lemma}
\begin{proof}
Set $\mu_{N,\varepsilon}:=\psi_\varepsilon*\mu_N$ and $\mu_\varepsilon:=\psi_\varepsilon*\mu$.
We first show that $\big(\mu_{N,\varepsilon}\big)_N$ converges narrowly to $\mu_\varepsilon.$ 
In particular, for any $\varphi\in C_b(\mathbb{R}^d)$, since $\psi_\varepsilon$ is a probability density, $ \|\Phi_\varepsilon\|_{L^\infty(\mathbb{R}^d)} \leq \|\varphi\|_{L^\infty(\mathbb{R}^d)}$ if we set 
\[
\Phi_\varepsilon(y) := \int_{\mathbb{R}^d} \varphi(y+z)\psi_\varepsilon(z)\,\ud z.
\]
Moreover, the continuity of $\varphi$ and the dominated convergence theorem imply that $\Phi_\varepsilon$ is continuous. Thus, $\Phi_\varepsilon\in C_b(\mathbb{R}^d)$. By Fubini's theorem and the narrow convergence of $\mu_N$,
\begin{equation*} 
 \int_{\mathbb{R}^d}\varphi(x)\,\ud\mu_{N,\varepsilon}(x) 
 = \int_{\mathbb{R}^d}  \Phi_\varepsilon(y)\,\ud\mu_N(y)\to \int_{\mathbb{R}^d} \Phi_\varepsilon(y)\,\ud\mu(y) =
  \int_{\mathbb{R}^d} \varphi(x)\,\ud\mu_\varepsilon(x),
\end{equation*}
which proves the narrow convergence of $\big(\mu_{N,\varepsilon}\big)_N$ to $\mu_\varepsilon$.

Next, since $\psi$ has finite second moment,
inequality~\eqref{4-1-1} gives
\[
M_2(\mu_\varepsilon) \leq 2M_2(\mu)+2\varepsilon^2M_2(\psi) <\infty, \qquad
\sup_N M_2(\mu_{N,\varepsilon})\leq 2\sup_N M_2(\mu_N)+2\varepsilon^2M_2(\psi) <\infty.
\]
Hence, the convolved measures satisfy the hypotheses of Lemma~\ref{lemma-7-stability}. Hence, 
\begin{equation*}
        \liminf_{N\to\infty}
        \mathcal{F}_{\varepsilon,U}(\mu_N)
        =
        \liminf_{N\to\infty}
        \mathcal{F}_U(\mu_{N,\varepsilon})
        \geq
        \mathcal{F}_U(\mu_\varepsilon)
        =
        \mathcal{F}_{\varepsilon,U}(\mu),
\end{equation*}
which completes the proof.
\end{proof}

\begin{remark}
We proved a more general result of lower semi-continuity compared to \cite{carrillo2019blob}. For the blob regularization of the internal energy, Carrillo, Craig, and Patacchini \cite[Proposition~3.9]{carrillo2019blob} establish lower semicontinuity with respect to narrow convergence for $m>1$, without assuming uniformly bounded second moments. For $m=1$, they prove lower semicontinuity for Gaussian kernels under $W_2$-convergence in $P_2(\mathbb{R}^d)$.
\end{remark}
\subsection{Convergence in large sample limit: $N\to\infty$, with $\varepsilon>0$ fixed}
In this part, we prove the convergence of particle schemes for blob regularized energy functional. We first establish the following convergence of confinement potential and interaction potential.

\begin{proposition}\label{lemma-17}For any $\mu \in P_{2}(\mathbb{R}^{d})$, let $\mu_{N}$ be its canonical $N$-particle approximation, under Assumption~\ref{assumption}, we have
\begin{equation*}
    \limsup_{N \to \infty} \mathcal{F}_{V}(\mu_{N}) \le \mathcal{F}_{V}(\mu),\quad \limsup_{N \to \infty} \mathcal{F}_{W}(\mu_{N}) \le \mathcal{F}_{W}(\mu).
\end{equation*}
\end{proposition}

\begin{proof}
For the energy with confinement potential $\mathcal{F}_{V}$, we split it into 
\begin{equation*}
    \mathcal{F}_{V}(\mu_{N}) = \int_{\mathbb{R}^{d}} V_{\mathrm{pos}} \ud \mu_{N} + \int_{\mathbb{R}^{d}} V_{\mathrm{res}} \ud \mu_{N}.
\end{equation*}
Since $V_{\mathrm{res}}$ is sub-quadratic, and $\mu_{N} \to \mu$ in $P_{2}(\mathbb{R}^{d})$, by part 3 of Lemma~\ref{lemma-3-p2}, we have
\begin{equation*}
    \lim_{N \to \infty} \int_{\mathbb{R}^{d}} V_{\mathrm{res}} \ud \mu_{N} = \int_{\mathbb{R}^{d}} V_{\mathrm{res}} \ud \mu.
\end{equation*}
For the convex term $V_{\mathrm{pos}}$, since $\mu_{N} := X_N \# \nu_{[0, 1]^{d}}$ where $X_N$ is the piece-wise constant approximation of the Brenier map $\nabla \phi_{\to \mu}$ given in the proof of Lemma~\ref{lemma-7-4}, we have
\begin{equation*}
    \int_{\mathbb{R}^{d}} V_{\mathrm{pos}}(x) \ud \mu_{N} = \int_{[0, 1]^{d}} V_{\mathrm{pos}}(X_N) \ud \nu_{[0, 1]^{d}}.
\end{equation*}
Note further that $\nabla \phi_{\rightarrow \mu } (\theta)= X (\theta) = x_\ell$ on $E_\ell$, by Jensen's inequality and the convexity of $V_{\mathrm{pos}}$, we have
\begin{equation*}
    \begin{aligned}
        \int_{[0, 1]^{d}} V_{\mathrm{pos}}(X_N) \ud \nu_{[0, 1]^{d}} &= \sum_{\ell = 1}^{N} \nu_{[0, 1]^{d}}(E_{\ell}) \times V_{\mathrm{pos}}(x_{\ell})\\
        &\le \sum_{\ell = 1}^{N} \nu_{[0, 1]^{d}}(E_{\ell}) \times \frac{1}{\nu_{[0, 1]^{d}}(E_{\ell})} \int_{E_{\ell}} V_{\mathrm{pos}}(\nabla \phi_{\to \mu}(\theta)) \ud \nu_{[0, 1]^{d}}\\
        &= \int_{[0, 1]^{d}} V_{\mathrm{pos}}(\nabla \phi_{\to \mu}(\theta)) \ud \nu_{[0, 1]^{d}} = \int_{\mathbb{R}^{d}} V_{\mathrm{pos}} \ud \mu
    \end{aligned}
\end{equation*}
Hence, $\limsup_{N \to \infty} \mathcal{F}_{V}(\mu_{N}) \le \mathcal{F}_{V}(\mu)$.

Likewise, for the interaction energy $\mathcal{F}_{W}$, we split it into
\begin{equation*}
    \mathcal{F}_{W}(\mu_{N}) = \frac{1}{2}\int_{\mathbb{R}^{d} \times \mathbb{R}^{d}} W_{\mathrm{pos}}(x - y) \ud \mu_{N}(x) \otimes \mu_{N}(y) + \frac{1}{2}\int_{\mathbb{R}^{d} \times \mathbb{R}^{d}} W_{\mathrm{res}}(x - y) \ud \mu_{N}(x) \otimes \mu_{N}(y).
\end{equation*}
Since $W_{\mathrm{res}}$ is sub-quadratic, we have,
\begin{equation*}
    W_{\mathrm{res}}(x - y) \le C(1 + |x - y|^{\beta}) \le C(1 + 2|x|^{\beta} + 2|y|^{\beta}),
\end{equation*}
and note that $\mu_{N} \otimes \mu_{N} \to \mu \otimes \mu$ in $P_{2}(\mathbb{R}^{d} \times \mathbb{R}^{d})$, we have
\begin{equation*}
    \lim_{N \to \infty} \int_{\mathbb{R}^{d} \times \mathbb{R}^{d}} W_{\mathrm{res}}(x - y) \ud \mu_{N}(x) \otimes \mu_{N}(y) = \int_{\mathbb{R}^{d} \times \mathbb{R}^{d}} W_{\mathrm{res}}(x - y) \ud \mu(x) \otimes \mu(y).
\end{equation*}
For the convex part $W_{\mathrm{pos}}$, since $\mu_{N} := X \# \nu_{[0, 1]^{d}}$, where $X$ is again the piecewise constant approximation of the Brenier map $\nabla \phi_{\to \mu}$ as given in the proof of Lemma~\ref{lemma-7-4}, we have
\begin{equation*}
    \int_{\mathbb{R}^{d} \times \mathbb{R}^{d}} W_{\mathrm{pos}}(x - y) \ud \mu_{N}(x) \otimes \mu_{N}(y) = \int_{[0, 1]^{d}} \int_{[0, 1]^{d}} W_{\mathrm{pos}}(X(\theta) - X(\underline{\theta})) \ud \theta \ud \underline{\theta}.
\end{equation*}
Using Jensen's inequality and the convexity of $W_{\mathrm{pos}}$, we get
\begin{equation*}
    \begin{aligned}
        &\int_{[0, 1]^{d}} \int_{[0, 1]^{d}} W_{\mathrm{pos}}(X(\theta) - X(\underline{\theta})) \ud \theta \ud \underline{\theta} = \sum_{\ell_{0} = 1}^{N} \sum_{\ell_{1} = 1}^{N} \nu(E_{\ell_{0}}) \nu(E_{\ell_{1}}) W_{\mathrm{pos}}(x_{\ell_{0}} - x_{\ell_{1}})\\
        &\le \sum_{\ell_{0} = 1}^{N} \sum_{\ell_{1} = 1}^{N} \nu(E_{\ell_{0}}) \nu(E_{\ell_{1}}) \frac{1}{\nu(E_{\ell_{0}})} \int_{E_{\ell_{0}}} W_{\mathrm{pos}}(\nabla \phi_{\to \mu}(\theta) - x_{\ell_{1}}) \ud \theta\\
        &\le \sum_{\ell_{0} = 1}^{N} \sum_{\ell_{1} = 1}^{N} \nu(E_{\ell_{0}}) \nu(E_{\ell_{1}}) \frac{1}{\nu(E_{\ell_{0}})} \frac{1}{\nu(E_{\ell_{1}})} \int_{E_{\ell_{0}}}\int_{E_{\ell_{1}}} W_{\mathrm{pos}}(\nabla \phi_{\to \mu}(\theta) - \nabla \phi_{\to \mu}(\underline{\theta})) \ud \theta \ud \underline{\theta}\\
        &= \int_{\mathbb{R}^{d} \times \mathbb{R}^{d}} W_{\mathrm{pos}}(x - y) \ud \mu(x) \otimes \mu(y).
    \end{aligned}
\end{equation*}
Consequently, $\limsup_{N \to \infty} \mathcal{F}_{W}(\mu_{N}) \le \mathcal{F}_{W}(\mu)$.
\end{proof}
We next establish continuity of the regularized internal energy under narrow convergence with uniformly bounded second moments. The argument is motivated by the analysis of generalized geodesics in Wasserstein space and the differentiation of energy functionals along such curves in \cite[Proposition~3.10]{carrillo2019blob}.

 
\begin{proposition}\label{theorem-7-diffusion} Let $m>d/(d+2)$, and let $\psi_\varepsilon$ be  mollifiers where $\psi\in C_b(\mathbb R^d)$ is a probability density of finite second moment as in Theorem \ref{theorem:blob}.
If $(\mu_n)_n \subset  P_2(\mathbb R^d)$ converges narrowly to $\mu \in  P_2(\mathbb R^d)$ and $(M_2(\mu_n))_n$ is uniformly bounded, then 
\[
\lim_{n\to\infty}\mathcal F_{\varepsilon,U}(\mu_n)=\mathcal F_{\varepsilon,U}(\mu).
\]
In particular, $\mathcal F_{\varepsilon,U}$ is continuous with respect to $W_2$-convergence.
\end{proposition}
\proof Let
\[
\rho_n=\psi_\varepsilon*\mu_n,\qquad\rho=\psi_\varepsilon*\mu.
\]
Because $\psi_\varepsilon$ is bounded and continuous, we have
\begin{equation}\label{eq:nov10.2026.0}
\rho_n(x)\rightarrow\rho(x)
\quad\text{for every }x\in\mathbb R^d.
\end{equation}
Taking also into account the second-moment estimate \eqref{4-1-1}, we have 
\[
0\leq \rho_n(x),\rho(x)\leq \|\psi_\varepsilon\|_\infty, \qquad \sup_n\int_{\mathbb R^d}|x|^2\rho_n(x) \ud x\leq A<\infty.
\]
For every $\alpha$ strictly between $d/(d+2)$ and $1$, H\"older's inequality gives, for $R\geq1$,
\[
\begin{aligned}
\int_{|x|>R}\rho_n(x)^\alpha \ud x&\leq\left(\int_{|x|>R}|x|^2\rho_n(x) \ud x\right)^\alpha \left(\int_{|x|>R}|x|^{-\frac{2\alpha}{1-\alpha}} \ud x\right)^{1-\alpha} \leq C_{d,\alpha}A^\alpha R^{d(1-\alpha)-2\alpha}.
\end{aligned}
\]
for uniform constants $C_{d,\alpha}$ and $A$. Since $\alpha>d/(d+2)$, we have
\begin{equation}\label{eq:nov10.2026.1}
\sup_n \int_{|x|>R} \rho_n(x)^a \ud x \rightarrow 0,\quad \text{ as } R\to\infty.
\end{equation}
In the following, to prove $\mathcal{F}_{U}(\rho_n) \to \mathcal{F}_{U}(\rho)$ (or equivalently, the convergence $\mathcal F_{\varepsilon,U}(\mu_n)\rightarrow\mathcal F_{\varepsilon,U}(\mu)$), we invoke Vitali's dominant convergence theorem. In particular, using the point-wise convergence and $L^{\infty}$ bound, we have
\begin{equation*}
    \int_{B_R}U_m(\rho_n) \ud x\rightarrow \int_{B_R}U_m(\rho) \ud x, \quad \text{ for any fixed } R>0
\end{equation*}
In the following, we distinguish the three cases to control the tail, based on the decay of density~\eqref{eq:nov10.2026.1}.

\hfill\break
(Case $d/(d+2)<m<1$). Recall $U_m(s)=\frac{1}{m-1}s^m$. Taking $\alpha = m$ in \eqref{eq:nov10.2026.1} and using the point-wise convergence in~\eqref{eq:nov10.2026.0}, we obtain
\[
\sup_n\int_{|x|>R}|U_m(\rho_n(x))| \ud x\rightarrow0.
\]
which proves $\mathcal F_{\varepsilon,U}(\mu_n)\rightarrow\mathcal F_{\varepsilon,U}(\mu)$ by combining the convergence on $B_R$ above.

\hfill\break
(Case $m=1$). Choose an exponent $\alpha \in(d/(d+2),1)$. For $0\leq s\leq\|\psi_\varepsilon\|_\infty$, we have (see for example, \citep[Page 6]{JKO}),
\[
|s\log s|\leq C(s^\alpha + s).
\]
Hence, by using~\eqref{eq:nov10.2026.1}, together with the following inequality from the second-momentum bound, 
\[
\int_{|x|>R}\rho_n(x) \ud x\leq \frac A{R^2},
\]
we can control the tails, and complete the proof based on the convergence on bounded balls.

\hfill\break
(Case $m>1$). Since the densities are uniformly bounded, we have
\[
|U_m(\rho_n)| \leq C\rho_n.
\]
for a uniform constant $C$. We invoke the second-moment bound $M_{2}(\rho_{n}) \le A$ to controls the tail, and get the convergence based on  the convergence on bounded balls.
\endproof


Combining Propositions~\ref{lemma-17} and~\ref{theorem-7-diffusion}, we get the following.
\begin{theorem}[$\Gamma$-convergence]Given $(\nu_{N})_{N = 1}^{\infty}$, which converges to $\nu$ in $P_{2}(\mathbb{R})$, we have the following $\Gamma$-convergence under the narrow convergence
\begin{equation*}
    \Gamma\text{-}\lim_{N \to \infty} \JKO_{\varepsilon, N}(\cdot, \nu_{N}) = \JKO_{\varepsilon}(\cdot, \nu)
\end{equation*}
\end{theorem}

\begin{proof}The $\liminf$ inequality is the direct result of stability lemmas in Lemmas~\ref{lemma-7-stability} and~\ref{lemma-13-stability}.

On the other hand, for any $\mu \in P_{2}(\mathbb{R}^{d})$, let $\mu_{N}$ be the canonical $N$-particle sampling of $\mu$ in $P_{2}(\mathbb{R}^{d})$ (in Definition~\ref{def-7}), we have, from Propositions~\ref{lemma-17} and~\ref{theorem-7-diffusion}:
\begin{equation*}
    \limsup_{N \to \infty} \mathcal{F}_{V}(\mu_{N}) \le \mathcal{F}_{V}(\mu),\quad \limsup_{N \to \infty} \mathcal{F}_{W}(\mu_{N}) \le \mathcal{F}_{W}(\mu),\quad \lim_{N \to \infty} \mathcal{F}_{\varepsilon, U}(\mu_{N}) = \mathcal{F}_{\varepsilon, U}(\mu)
\end{equation*}
Since $W_{2}(\mu, \mu_{\tk})$ are continuous in $P_{2}(\mathbb{R}^{d})$, we have
\begin{equation*}
    \limsup_{N \to \infty} \JKO_{\varepsilon, N}(\mu_{N}, \nu_{N}) \le \limsup_{N \to \infty} \JKO_{\varepsilon}(\mu, \nu)\,,
\end{equation*}
which concludes the proof.
\end{proof}

We now state the main theorem in this subsection. The proof is the direct result of the above $\Gamma$-convergence of the JKO functional, which can be the deduced the same way as that of the one-dimensional proof. 
\begin{theorem}[blob method]\label{theorem:blob} Suppose that $\psi\in C_b(\mathbb R^d)$ is a probability density of finite second moment. 
Given initial value $\mu_{0} \in P_{2}(\mathbb{R}^d)$ and let $\mu_{0, N} := \frac{1}{N}\sum_{\ell = 1}^{N} \delta[x_{0, \ell}] \in P_{2, N}(\mathbb{R}^d)$ be a sampling of $\mu_{0}$ such that $\mu_{0, N} \to \mu_{0}$ in $P_{2}(\mathbb{R}^{d})$, we have
\begin{equation*}
    \mu_{\tk, N} \to \mu_{\varepsilon, \tk} \text{ in } P_{2}(\mathbb{R}^d) \text{ as } N \to \infty \text{ for all } \tk \in \mathbb{N}.
\end{equation*}
Here, $\mu_{\tk, N}$ is the $N$-particle solution to \eqref{JKO-blob}, and $\mu_{\varepsilon, \tk}$ is the solution of the $\varepsilon$-regularized JKO scheme \eqref{JKO_eps}.
\end{theorem}

\subsection{Simultaneous particle and regularization limits: $N\to\infty$, $\varepsilon_N\to0$, under $N\varepsilon_N^d\to\infty$}
The goal of this subsection is to prove Theorem \ref{theorem:blob:bis}.

\proof Suppose
\[
\mu_N\in P_{2,N}(\mathbb R^d),\quad \mu_N\rightharpoonup\mu,\quad \sup_NM_2(\mu_N)<\infty, \quad 
\widehat\mu_N=(\psi_{\varepsilon_N}*\mu_N)\mathcal L^d, 
\]
so that one has (see for example \cite{GangboTudorascu2014} and also~\eqref{4-1-1} above)
\[
W_2^2(\widehat\mu_N,\mu_N)\le\varepsilon_N^2 \int_{\mathbb R^d}|z|^2\psi(z) \ud z, \quad M_2(\widehat\mu_N)\le 2M_2(\mu_N)+2\varepsilon_N^2M_2(\psi).
\]
Applying Lemma \ref{lemma-7-stability}, we have
\[
\liminf_{N\to\infty}\mathcal F_{U}(\widehat\mu_N)\ge\mathcal F_{U}(\mu).
\]
The same lemma applies to the confinement and interaction energies evaluated on $\mu_N$. Consequently, 
\[
\liminf_{N\to\infty}\mathcal F_N(\mu_N)\ge\mathcal F(\mu),
\]
which proves the $\Gamma$-liminf inequality.

In the following, we first construct recovery sequences for smooth compactly supported density, and then pass them to general measures via density argument.

{\bf Step 1.} (Smooth compactly supported measures $\mu$). Suppose 
\[
\mu=\rho \ud x, \quad \rho\in C_c^\infty(\mathbb R^d), \quad \mu_N=\frac1N\sum_{\ell=1}^N\delta_{X_\ell}.
\]
where $\rho\ge0$, $\int_{\mathbb R^d}\rho \ud x=1$ and $X_1,\ldots,X_N$ are $N$ independent random variables with law $\mu$.  Define
\[
\rho_{N,\varepsilon}=\psi_\varepsilon*\mu_N\quad \text{and} \quad \rho_\varepsilon=\psi_\varepsilon*\rho.
\]
Since the $X_\ell$ are independent, the centered random variables
\[
Y_\ell(x):=\psi_\varepsilon(x-X_\ell)-\rho_\varepsilon(x)
\]
are independent and have mean zero. Therefore, for every $x\in\mathbb R^d$, we have 
\begin{equation*}
\mathbb E\left[|\rho_{N,\varepsilon}(x)-\rho_\varepsilon(x)|^2\right]=\sum_{\ell=1}^N \frac{1}{N^2}\mathbb E\left[|Y_\ell(x)|^2\right]
= \frac{1}{N} \mathbb E\left[|Y_1(x)|^2\right].
\end{equation*}
Furthermore,
\begin{equation}\label{eq:nov12.2026.9}
\mathbb E\left[|Y_1(x)|^2\right]= \int_{\mathbb R^d}\psi_\varepsilon(x-y)^2 \ud \mu(y)-\rho_\varepsilon(x)^2 \leq \int_{\mathbb R^d}\psi_\varepsilon(x-y)^2 \ud \mu(y).
\end{equation}
Integrating~\eqref{eq:nov12.2026.9} over $x\in\mathbb R^d$ and applying Tonelli's theorem gives
\begin{equation}\label{eq:nov12.2026.9bis}
\int_{\mathbb R^d }\mathbb E\left[|Y_1(x)|^2\right] \ud x \leq \int_{\mathbb R^d } \int_{\mathbb R^d}\psi_{\varepsilon_N}(x-y)^2\,\rho(y) \ud y \ud x=\|\psi_{\varepsilon_N}\|_{L^2}^2={\|\psi\|_{L^2}^2 \over \varepsilon_N^d}.
\end{equation}
Consequently, we have
\begin{equation}\label{eq:nov12.2026.1}
\mathbb E\|\rho_{N,\varepsilon_N}-\rho_{\varepsilon_N}\|_{L^2}^2 \leq {1\over N \varepsilon_N^{d}}\|\psi\|_{L^2}^2.
\end{equation}
We can then conclude that 
\[
\mathbb E\|\rho_{N,\varepsilon_N}-\rho\|_2^2\leq\frac{\|\psi\|_{L^2}^2}{N\varepsilon_N^d}+\|\psi_{\varepsilon_N}*\rho-\rho\|_2^2 \rightarrow 0.
\]
Hence
\begin{equation}\label{eq:nov12.2026.2}
\rho_{N,\varepsilon_N} \rightarrow \rho \text{ in }L^2(\mathbb R^d)\quad \text{in probability}.
\end{equation}
Because $\rho$ and $\psi$ are compactly supported, all these densities are supported in a common compact set for sufficiently large $N.$ 

We claim that 
\begin{equation}\label{eq:nov12.2026.10b}
\int_{\mathbb R^d} U(\rho_{N,\varepsilon_N}) \ud x \rightarrow\int_{\mathbb R^d} U(\rho) \ud x\quad \text{ in probability}. 
\end{equation}
For the cases where $d/(d+2)<m\le2$, the entropy function $U$ has a growth rate less than quadratic. Hence, we can prove~\eqref{eq:nov12.2026.10b} using \eqref{eq:nov12.2026.2} through H\"older-type inequalities. For the rest cases where $m > 2$, we observe that 
\[
\mathbb E\Big[|\rho_{N,\varepsilon_N}(x)-\rho_{\varepsilon_N}(x)|^m\Big]= \mathbb E\bigg|\frac1N\sum_{\ell=1}^NY_\ell(x)\bigg|^m.
\]
Since the random variables $Y_\ell(x)$ are i.i.d., and centered, we can apply Rosenthal's inequality to find a constant $C_m$ depending only on $m$ such that  
\begin{equation}\label{eq:nov12.2026.10}
\mathbb E\Big[|\rho_{N,\varepsilon_N}(x)-\rho_{\varepsilon_N}(x)|^m\Big] \le C_m\bigg[N^{1-m}\mathbb E|Y_1(x)|^m+N^{-m/2}\bigl(\mathbb E|Y_1(x)|^2\bigr)^{m/2}\bigg].
\end{equation}
On the other hand, since $Y_\ell(x):=\psi_\varepsilon(x-X_\ell)-\rho_\varepsilon(x)$, we have 
\[
\mathbb E|Y_1(x)|^m \leq 2^{m-1} \big(\mathbb E|\psi_{\varepsilon_N}(x-X_1)|^m +|\rho_\varepsilon(x)|^m\big),
\]
By applying Jensen's inequality to $|\rho_\varepsilon(x)|^m$, we obtain   
\[
\mathbb E|Y_1(x)|^m \leq 2^{m}\ \mathbb E|\psi_{\varepsilon_N}(x-X_1)|^m = 
2^{m} \int_{\mathbb R^d} |\psi_{\varepsilon_N}(x-y)|^m\rho(y) \ud y.
\]
Integrating the above inequality over $x \in \mathbb R^d$, we have 
\begin{equation}\label{eq:nov12.2026.8}
\int_{\mathbb R^d} \mathbb E|Y_1(x)|^m dx\leq 2^m \|\psi_{\varepsilon_N}\|_{L^m}^m={2^{m} \over \varepsilon_N^{d(m-1)}} \|\psi\|_{L^m}^m.
\end{equation}

Let
\[
g_N(x):=\int_{\mathbb R^d}\psi_{\varepsilon_N}(x-y)^2 \ud \mu(y).
\]
Here we note that
\[
\|g_N\|_{L^1}=\|\psi_{\varepsilon_N}\|_{L^2}^2=\varepsilon_N^{-d}\|\psi\|_{L^2}^2, \quad 
\|g_N\|_{L^\infty}\le\|\rho\|_{L^\infty}\|\psi_{\varepsilon_N}\|_{L^2}^2=\|\rho\|_{L^\infty}\varepsilon_N^{-d}\|\psi\|_{L^2}^2.
\]
Hence, we can then use~\eqref{eq:nov12.2026.9} to infer that
\begin{equation}\label{eq:nov12.2026.8-1}
\int_{\mathbb R^d }\big(\mathbb E\left[|Y_1(x)|^2\right]\big)^{m \over 2} \ud x \leq \int_{\mathbb R^d }  g_N(x)^{m \over 2} \ud x \leq 
\|g_N\|_{L^\infty}^{{m \over 2}-1} \|g_N\|_{L^\infty} \leq \|\rho\|_{L^\infty}^{{m\over 2}-1}
\|\psi\|_{L^2}^{m}\varepsilon_N^{-{dm\over 2}}.
\end{equation}
Combining estimates~\eqref{eq:nov12.2026.8} and~\eqref{eq:nov12.2026.8-1} into inequality \eqref{eq:nov12.2026.10}, we can find a constant $C_{m,\rho,\psi}$ depending only on $m$, $\rho$ and $\psi$ such that 
\begin{equation}\label{eq:nov12.2026.7}
\mathbb E\left[\|\rho_{N,\varepsilon_N}-\rho_{\varepsilon_N}\|^m_{L^m}\right] \leq C_{m,\rho,\psi}\left[(N\varepsilon_N^d)^{-(m-1)}+(N\varepsilon_N^d)^{-m/2}\right]  \longrightarrow 0,
\end{equation}
when $N\varepsilon_N^d \rightarrow +\infty.$ Since $\rho_{\varepsilon_N}\to\rho$ in $L^m(\mathbb{R}^d)$, we obtain that $\mathbb E\|\rho_{N,\varepsilon_N}-\rho\|_m^m\to0.$ This is sufficient to conclude that \eqref{eq:nov12.2026.10b} holds.

Furthermore, because the common support of $(\mu_N)_N$ is compact, we have the sequence of empirical measures $(\mu_N)_N$ converges to $\mu$ in $ P_{2}(\mathbb{R}^d)$ almost surely. The confinement and interaction energies also converge, by continuity of $V,W$ on compact sets. In conclusion, 
we have proved that, for every smooth compactly supported probability density $\mu$, there exists a sequence $\mu_N\in P_{2,N}(\mathbb{R}^d)$ for all $N$, satisfying
\begin{equation}\label{eq:nov12.2026.11}
W_2(\mu_N,\mu)\to0,\quad \mathcal F_N(\mu_N)\to\mathcal F(\mu).
\end{equation}

{\bf Step 2.} (Extension to general finite-energy measures). Suppose $\mu\in P_2(\mathbb R^d)$ and $\mathcal F(\mu)<\infty$. Let
\[
\mu=\rho \ud x+\mu^s,\qquad a_R=\mu(B_R),\qquad\mu^R=\frac{\mu|_{B_R}}{a_R}.
\]
We note that $\mu^R\to\mu$ in $W_2$. For $m\ne1$, we have 
\[
\mathcal F_{U}(\mu^R) =\frac{a_R^{-m}}{m-1}\int_{B_R}\rho(x)^m \ud x.
\]
Since $a_R\to1$ and $\mathcal F_{U}(\mu)$ is finite, we have $\mathcal F_{U}(\mu^R)\to\mathcal F_{U}(\mu)$. We also draw similar conclusions for $m=1.$ Furthermore, and the potential energies $\mathcal{F}_{V}(\mu^R)$ and $\mathcal{F}_{W}(\mu^R)$ converge by integrability and the assumed growth structure.

Fix $R$ and let $\eta_\delta$ be a compactly supported smooth probability kernel, we use the mollification
\[
\mu^{R,\delta}=\eta_\delta*\mu^R.
\]
Convexity of $U$ and the recession function convention, give $\mathcal F_U(\mu^{R,\delta}) \le \mathcal F_U(\mu^R).$ 
Its lower semicontinuity and $W_2(\mu^{R,\delta},\mu^R)\to0$ then imply that $\mathcal F_U(\mu^{R,\delta})$ converges to $\mathcal F_U(\mu^R)$ as $\delta$ decreases to $0$.

For fixed $R$, since the supports of the measures remain in a common compact set when $\delta$ is small,  we obtain that the confinement and interaction energies converge.  A diagonal argument produces the desired smooth compactly supported recovery targets.

Let $\nu_j=\rho_j \ud x$, where $\rho_j\in C_c^\infty(\mathbb R^d)$ are non-negative probability densities such that $(\nu_j)_j$ converges to $\mu$ in $W_2$ and $(\mathcal F(\nu_j))_j$ converges to $\mathcal F(\mu)$ as $j \to +\infty.$
By using the result in \textbf{Step 1}, for each fixed $j$, there exist empirical measures $\mu_N^{(j)}\in P_{2,N}(\mathbb R^d)$ such that
\[
W_2(\mu_N^{(j)},\nu_j)\to0,\quad
\mathcal F_N(\mu_N^{(j)})\to\mathcal F(\nu_j),\quad \text{ as $N\to\infty$.}
\]
Choose a  sequence $(N_j)_j$ increasing to $+\infty$ such that, for every $N\ge N_j$,
\[
W_2(\mu_N^{(j)},\nu_j)+\left|\mathcal F_N(\mu_N^{(j)})-\mathcal F(\nu_j)\right|\le\frac1j.
\]
For $N_j\le N<N_{j+1}$, set $\mu_N=\mu_N^{(j)}$. Then 
\[
\begin{aligned}
W_2(\mu_N,\mu)\le\frac1j+W_2(\nu_j,\mu)\rightarrow0,\quad 
\left|\mathcal F_N(\mu_N)-\mathcal F(\mu)\right| \le\frac1j+\left|\mathcal F(\nu_j)-\mathcal F(\mu)\right| \rightarrow0.
\end{aligned}
\]
Consequently, the empirical measure $(\mu_N)$ is a $W_2$-convergent recovery sequence. Together with the $\Gamma$-liminf inequality, we have $\Gamma\text{-}\lim_{N \to \infty}\mathcal F_N = \mathcal F$.

{\bf Step 3.} (Convergence of the JKO minimizers). Fix $\tau>0$, and suppose that $\nu_N\to\nu$ in $W_2$. Set
\[
\begin{aligned}
\mathcal J_N(\mu)=\mathcal F_N(\mu)+\frac1{2\tau}W_2^2(\mu,\nu_N), \quad 
\mathcal J(\mu) =\mathcal F(\mu)+\frac1{2\tau}W_2^2(\mu,\nu).
\end{aligned}
\]
Since $\nu_N\to\nu$ in $W_2$, the $\Gamma$-convergence of $\mathcal F_N$ established above implies
\[
\Gamma\text{-}\lim_{N\to\infty}\mathcal J_N=\mathcal J.
\]
To verify the uniform coercivity required by Lemma \ref{lemma-12}, we note that for every $\mu\in P_{2,N}(\mathbb R^d)$,
\[
M_2((\psi_{\varepsilon_N}*\mu) \ud x) \le 2M_2(\mu)+2\varepsilon_N^2M_2(\psi).
\]
By Lemmas \ref{lemma-5-internal} and \ref{lemma-6-well}, and because $\varepsilon_N\le1$ for all sufficiently large $N$, we can find constants $C>0$ and $0<\beta<2$, independent of $N$, such that
\[
\mathcal F_N(\mu)\ge-C\bigl(1+M_2(\mu)^{\beta/2}\bigr).
\]
Consequently,
\begin{equation}\label{w-c-5-1}
\begin{aligned}
\mathcal J_N(\mu)\ge \frac1{4\tau}M_2(\mu)-\frac1{2\tau}M_2(\nu_N)-C\bigl(1+M_2(\mu)^{\beta/2}\bigr).
\end{aligned}
\end{equation}
Since $\sup_NM_2(\nu_N)<\infty$ and $\beta<2$, the sublevel sets of $\mathcal J_N$ have uniformly bounded second moments. For each fixed $N$, the functional $\mathcal J_N$ is lower semi-continuous with respect to the particle positions (by Lemma~\ref{lemma-7-stability}). Its sublevel sets are bounded by the preceding coercivity estimate in~\eqref{w-c-5-1}. Therefore, the minimization problem in $(\mathbb R^d)^N$ gives the existence of a minimizer. Let $\mu_N^*$ be a minimizer of $\mathcal J_N$, and suppose the minimizer $\mu^*$ of $\mathcal J$ is unique. By Lemma  \ref{lemma-12}, $W_2(\mu_N^*,\mu^*)\rightarrow0.$ 

Finally, assume that the sequence $(\mu_{0,N})_N$ of initial empirical measures converges to $\mu_0$ in $ P_2(\mathbb R^d)$ as $N$ tends to $+\infty.$  Applying the preceding conclusion to $\nu_N=\mu_{\tk-1,N}$ and $\nu=\mu_{\tk-1}$, and arguing inductively, we obtain that for every fixed $\tk\geq 1,$ $(\mu_{\tk,N})_N$ converges to $\mu_\tk$ in $P_2(\mathbb R^d)$ as $N$ tends to $+\infty$. This completes the proof.\endproof

\section{Transport map approach}\label{sec:5}
This section is devoted to establishing the convergence in Theorem~\ref{theorem-3-neural}. At each step, the solution is defined by minimizing the following discrete JKO problem.
\begin{definition}[discrete JKO functional, transport map]\rm Given $N \in \mathbb{N}$, we define
\begin{equation} \label{1002}
    \JKO_{N}(\phi; \nu^{\mathrm{ac}}_{N}, \nu_{N}) := \mathcal{F}_{N}(\phi; \nu^{\mathrm{ac}}_{N}, \nu_{N}) + \frac{1}{2\tau}W_{2}\big((\nabla \phi) \# \nu_{N}, \nu_{N}\big)^{2}
\end{equation}
where
\begin{equation*}
    \mathcal{F}_{N}(\phi; \nu^{\mathrm{ac}}_{N}, \nu_{N}) := \mathcal{F}_{U, N}(\phi; \nu^{\mathrm{ac}}_{N}, \nu_{N}) + \mathcal{F}_{V}\big(\nabla \phi) \# \nu_{N}\big) + \mathcal{F}_{W}\big((\nabla \phi) \# \nu_{N}\big),
\end{equation*}
and the discrete internal energy is given by
\begin{equation*}
    \mathcal{F}_{U, N}(\phi; \nu^{\mathrm{ac}}_{N}, \nu_{N}) := \int_{\mathbb{R}^{d}} S\bigg(\frac{\nu^{\mathrm{ac}}_{N}(x)}{\det(\nabla^{2} \phi(x))}\bigg) \ud \nu_{N}\,, ~\text{ where $S(s) = U(s)/s$.}
\end{equation*}
\end{definition}

As we mentioned in Assumption \ref{assump3}, we consider initial values $\mu_{0} \in C\big(\mathbb{R}^{d}\big)$ and satisfies \eqref{lower-bound-3}. By the explicit expression of $S(s)=U(s)/s$ and Assumption 3, there exists $C>0$ such that
\begin{equation}\label{5-1-bound}
-C(1+|x|^2)\le S(\mu_0(x))\le C, \qquad x\in\mathbb R^d.
\end{equation}

We minimize $\JKO_N(\phi;\nu_N^{\mathrm{ac}},\nu_N)$ over the admissible class $C_{\Lambda}^{2,\alpha}(\mathbb{R}^d)$ defined in~\eqref{Calpha}, where $\Lambda$, $\alpha$, and the constants $C_R$ are fixed independently of the potential $\phi$. The admissible class is closed under locally uniform convergence. In particular, if $\sup_{N} [\nabla^2 \phi_{N}]_{\alpha, B_{R}(0)}$ is bounded for all $R$, then $\phi_N \to \phi$ locally uniformly implies
\begin{equation*}
    \phi \in C^{2, \alpha}_{\mathrm{loc}}(\mathbb{R}^d),\quad \text{ and } \phi_N \to \phi \text{ in } C_{\mathrm{loc}}^{2, \beta}(\mathbb{R}^d) \text{ for } 0 < \beta < \alpha,
\end{equation*}
based on interpolation inequalities between $L^{\infty}$ norm and H\"older norms (see, for example, \cite[Lemmas 6.35-6.36]{gilbarg1977elliptic}). The $C_{\mathrm{loc}}^{2, \beta}$ convergence ensures that $\nabla^2\phi_N$ converges locally uniformly, and therefore, the pointwise Hessian bounds pass to the limit:
\[
    \Lambda I_d
    \leq \nabla^2\phi(x)
    \leq \Lambda^{-1}I_d
    \quad\text{for all }x\in\mathbb{R}^d,
\]
which implies that $\phi \in C_{\Lambda}^{2, \alpha}(\mathbb{R}^d)$.

\subsection{Existence of minimizer}
\begin{theorem}\label{theorem-1002-existence}For each $N$, and given $\nu_{N} \in P_{2}(\mathbb{R}^d)$ and $\nu_{N}^{\mathrm{ac}} \in P_{2}(\mathbb{R}^d) \cap C(\mathbb{R}^d)$ with
\begin{equation*}
    -C(1 + |x|^{2}) \le S(\nu^{\mathrm{ac}}_{N}) \leq C, 
\end{equation*}
for a positive constant $C$. Then there exists a minimizer $\phi_{N}$ to the discrete JKO problem~\eqref{1002}.
\end{theorem}

We first estimate the second moment for the sub-level sets associated with the discrete JKO problem \eqref{1002}, and establish a lemma in analogous with Lemma~\ref{lemma-11-coercive}, that these sublevel sets are contained in a class of measures with uniformly bounded second moments.

We leverage the ellipticity:
\begin{equation} \label{ellip}
    \Lambda I_{d} \le \nabla^2 \phi \le \Lambda^{-1}I_d.
\end{equation}
In particular, this condition implies that, for any $\phi \in C_{\Lambda}^{2, \alpha}(\mathbb{R}^{d})$, one has $\phi(x) \to \infty$ as $|x| \to \infty$ and $\phi$ is strongly convex. Consequently, for each $\phi$ in $C_{\Lambda}^{2, \alpha}(\mathbb{R}^{d})$, it admits a unique minimizer. \vskip0.05in

\centerline{We denote the unique minimizer of $\phi$ by $x_{0}(\phi) \in \mathbb{R}^{d}$.}
\begin{lemma}[Bounds on $x_{0}(\phi)$]\label{lemma:15-equi} For any $\mu \in P_{2}(\mathbb{R}^{d}),$ and $\phi \in C_{\Lambda}^{2, \alpha}(\mathbb{R}^{d})$, 
\begin{equation*}
    \tfrac{1}{2}\Lambda^{2} |x_{0}(\phi)|^{2} - 2M_2(\mu)\Lambda^{-2} \le M_{2}((\nabla \phi)\# \mu) \le 2\Lambda^{-2} |x_{0}(\phi)|^{2} + 2M_2(\mu)\Lambda^{-2}.
\end{equation*}
\end{lemma}
\begin{proof}
Let $x_0:=x_0(\phi)$. Since $\nabla\phi(x_0)=0$, the fundamental theorem of calculus on polar coordinate $r:=|x-x_0|$ and $v:=(x-x_0)/r$ gives
\begin{equation*}
\nabla\phi(x)
=\int_0^r \nabla^2\phi(x_0+tv)v\ud t.
\end{equation*}
Using the ellipticity \eqref{ellip}
and $|v|=1$, we obtain 
\[
\Lambda r\le \int_0^r\langle v,\nabla^2\phi(x_0+tv)v\rangle\ud t=\langle v,\nabla\phi(x)\rangle= \int_0^r|\nabla^2\phi(x_0+tv)v|\ud t
\le\Lambda^{-1}r.
\]
Maximizing or minimizing over $x$ such that $r:=|x-x_0|$ we have 
\begin{equation*}
\Lambda|x-x_0(\phi)|
\le|\nabla\phi(x)|
\le\Lambda^{-1}|x-x_0(\phi)|.
\end{equation*}
Since $M_{2}((\nabla \phi)\# \mu) = \int_{\mathbb{R}^{d}} |\nabla \phi|^{2} \ud \mu$, then
\begin{equation*}
    \int_{\mathbb{R}^{d}} |\nabla \phi|^{2} \ud \mu \le \Lambda^{-2} \int_{\mathbb{R}^{d}} |x - x_{0}(\phi)|^{2} \ud \mu \le 2\Lambda^{-2} \int_{\mathbb{R}^{d}} (|x|^{2} + |x_{0}(\phi)|^{2}) \ud \mu.
\end{equation*}
Likewise, we obtain the necessary lower bound and conclude the proof.
\end{proof}
\begin{proposition}[Concentration estimate for the discrete JKO] \label{theorem-10-concentration} Suppose that $\nu_{N}, \nu^{\mathrm{ac}}_{N} \in P_{2}(\mathbb{R}^{d})$ and 
\begin{equation}\label{9-1-bound}
    -C(1 + |x|^{2}) \le S(\nu^{\mathrm{ac}}_{N}) \leq C
\end{equation}
for some constant $C>0$.  For any  $\phi \in C_{\Lambda}^{2}(\mathbb{R}^{d})$ that satisfies  $\JKO_{N}(\phi; \nu^{\mathrm{ac}}_{N}, \nu_{N}) \le E$ for some constant $E >0$, there exists a constant $R = R(\nu^{\mathrm{ac}}_{N}, \nu_{N}, E, \Lambda)$ such that $ |x_{0}(\phi)| \le R.$
\end{proposition}
\begin{proof}
Denote $\mu_{N} = (\nabla \phi) \# \nu_{N}$. For $W_{2}$, $V$, and $W$, we have
\begin{equation*}
    W_{2}(\mu_{N}, \nu_{N})^{2} \ge \tfrac{1}{2}M_{2}(\mu_{N}) - M_{2}(\nu_{N}).
\end{equation*}
\begin{equation*}
    \mathcal{F}_{V}(\mu_{N}) \ge -C(1 + M_{2}(\mu_{N})^{\beta/2}),\quad \mathcal{F}_{W}(\mu_{N}) \ge -2C(1 + M_{2}(\mu_{N})^{\beta/2}),
\end{equation*}
where the inequality of $W_{2}$ is given by \citep[Equation (17)]{JKO}, and the inequalities for confinement potential and interaction potential follow from Lemma~\ref{lemma-6-well}.

For the internal energy, we have 
\begin{equation*}
    \begin{aligned}
        \mathcal{F}_{U, N}(\phi; \nu^{\mathrm{ac}}_{N}, \nu_{N}) &= \int_{\mathbb{R}^{d}} S\bigg(\frac{\nu^{\mathrm{ac}}_{N}(x)}{\det(\nabla^{2} \phi(x))}\bigg) \ud \nu_{N}\\
        &= \int_{\mathbb{R}^{d}} S\bigg(\frac{\nu^{\mathrm{ac}}_{N}(x)}{\det(\nabla^{2} \phi(x))}\bigg) \ud \nu^{\mathrm{ac}}_{N} + \int_{\mathbb{R}^{d}} S\bigg(\frac{\nu^{\mathrm{ac}}_{N}(x)}{\det(\nabla^{2} \phi(x))}\bigg) \ud \big(\nu_{N} - \nu^{\mathrm{ac}}_{N}\big).
    \end{aligned}
\end{equation*}
Since $\Lambda^{d} \le \det(\nabla^{2} \phi(x)) \le \Lambda^{-d}$ and using~\eqref{9-1-bound}, we have 
\begin{equation*}
    \bigg|S\bigg(\frac{\nu^{\mathrm{ac}}_{N}(x)}{\det(\nabla^{2} \phi(x))}\bigg)\bigg| \le C(1 + |x|^{2})
\end{equation*}
for some constant $C$. As a result, we can deduce that the  error of finite sampling
\begin{equation*}
    e_{N} := \int_{\mathbb{R}^{d}} S\bigg(\frac{\nu^{\mathrm{ac}}_{N}(x)}{\det(\nabla^{2} \phi(x))}\bigg) \ud (\nu_{N} - \nu^{\mathrm{ac}}_{N})
\end{equation*}
is bounded below by a constant $C(\nu_N, \nu_N^\text{ac})$, which only depends on $M_{2}(\nu_{N})$ and $M_{2}(\nu^{\mathrm{ac}}_{N})$. Therefore, we have
\begin{equation*}
    \begin{aligned}
        \mathcal{F}_{U, N}(\phi; \nu^{\mathrm{ac}}_{N}, \nu_{N}) &= e_{N} + \mathcal{F}_{U}((\nabla \phi) \# \nu^{\mathrm{ac}}_{N}) \\
        &\ge C(\nu_{N}, \nu^{\mathrm{ac}}_{N}) + \begin{cases}
        C(1, d) - \frac{d}{2}\log M_{2}((\nabla \phi) \# \nu^{\mathrm{ac}}_{N}) &\text{ when } m = 1\,;\\
        C(m, d)M_{2}((\nabla \phi) \# \nu^{\mathrm{ac}}_{N})^{\frac{d}{2}(1 - m)} &\text{ when } m \neq 1\,.
    \end{cases}
    \end{aligned}
\end{equation*}
Using Lemma~\ref{lemma:15-equi} to express all the $M_{2}((\nabla \phi) \# \nu^{\mathrm{ac}})$ and $M_{2}((\nabla \phi) \# \nu_{N})$ in terms of $|x_{0}(\phi)|$, we have
\begin{equation*}
    \JKO_{N}(\phi; \nu^{\mathrm{ac}}_{N}, \nu_{N}) \ge \tfrac{1}{8\tau} \Lambda^{2} |x_{0}(\phi)|^{2} - C(1 + |x_{0}(\phi)|^{\beta})
\end{equation*}
for $\beta < 2$, which concludes the proof.
\end{proof}
\begin{lemma}[Compactness of $C_{\Lambda}^{2, \alpha}(\mathbb{R}^{d})$]\label{lemma-18-compact} Let $(\phi_{N})_{N = 1}^{\infty}$ be a sequence of functions in $C_{\Lambda}^{2, \alpha}(\mathbb{R}^{d})$ such that $x_{0}(\phi_{N})$ is uniformly bounded. Then up-to the extraction of subsequence, there exists a function $\phi \in C_{\Lambda}^{2, \alpha}(\mathbb{R}^{d})$ such that $\nabla \phi_{N} \to \nabla \phi$ locally uniformly, and in addition, $\nabla^{2}\phi_{N} \to \nabla^{2}\phi$ in $C^{0, \beta}(B_{R}(0))$ for all $R$ with $0 < \beta < \alpha$.
\end{lemma}

\begin{proof} Let $(\phi_{N})_{N = 1}^{\infty}$ be a sequence of convex function $C_{\Lambda}^{2, \alpha}(\mathbb{R}^{d})$ such that $x_{0}(\phi_{N})$ is uniformly bounded. Since $\phi$ and $\phi + c$ has the same gradient and Hessian, we consider the following normalization, without loss of generality,
\begin{equation*}
    \phi_{N}\big(x_{0}(\phi_{N})\big) = 0.
\end{equation*}
Let $R_{0} = \sup_{N} x_{0}(\phi_{N})$. Given any $R$, one has
\begin{equation*}
    |\phi_{N}(x)| \le \frac{|x - x_{0}(\phi_{N})|^{2}}{2\Lambda} \le \frac{R^{2} + R_{0}^{2}}{\Lambda},\quad |\nabla \phi_{N}(x)| \le \frac{|x - x_{0}(\phi_{N})|}{\Lambda} \le \frac{R + R_{0}}{\Lambda}.
\end{equation*}
for any $|x| \le R$. By Arzelà-Ascoli theorem, for each $R < \infty$, $\big(\phi_{N}\big)_N$ converges uniformly to $\phi$ on $B_{R}(0)$, 
up to the extraction of a subsequence. Now, to extend the convergence to the domain $\mathbb{R}$, we consider a sequence of balls $B_{j}(0)$ for $j = 1, \cdots, \infty$. By a standard diagonalization argument, $\big(\phi_{N}\big)_N$ converges uniformly to $\phi$ on $B_{j}(0)$  for all  $j \in \mathbb N.$ Extraction of a subsequence, we conclude that $\big(\phi_{N}\big)_N$ converges locally uniformly to $\phi.$ 
Based on interpolation inequality in H\"older spaces, the $C^{2, \alpha}_{\mathrm{loc}}(\mathbb{R}^{d})$ H\"older regularity together with the local uniform convergence, implies that $\phi_{N} \to \phi$ in $C^{2, \beta}_{\mathrm{loc}}(\mathbb{R}^{d})$ for $0 < \beta < \alpha$.
\end{proof}

\begin{proof}[Proof of Theorem~\ref{theorem-1002-existence}] According to Proposition~\ref{theorem-10-concentration}, the JKO functional is bounded from below. We then let $(\phi_{j})_j$ be a minimizing sequence for $\JKO_{N}(\cdot, \nu_{N}^{\mathrm{ac}},\nu_N)$.

Due to Proposition~\ref{theorem-10-concentration}, $(x_0(\phi_{j}))_j$ is uniformly bounded. We fix $\phi_j(0) = 0$ for all $j$ without loss of generality, because the JKO functional is invariant under $\phi \mapsto \phi + c$ with any constant $c$. According to Lemma~\ref{lemma-18-compact}, there exists a limit $\phi_\ast \in C_{\Lambda}^{2, \alpha}(\mathbb{R}^d)$ such that
\begin{equation*}
    \phi_{j} \to \phi_\ast \text{ in } C^{2, \beta}_{\mathrm{loc}}(\mathbb{R}^d) \text{ as } j \to \infty
\end{equation*}
up to the extraction of a subsequence. By the local uniform convergence of $\nabla^2 \phi_j \to \nabla^2 \phi_\ast$, the ellipticity $\Lambda^{d} \le \det(\nabla^2 \phi_j) \le \Lambda^{-d}$, and the continuity of $S$ in $(0, \infty)$, we have, using dominant convergence:
\begin{equation*}
    \lim_{j \to \infty} \mathcal{F}_{U, N}(\phi_j; \nu^{\mathrm{ac}}_{N}, \nu_{N}) = \lim_{j \to \infty} \int_{\mathbb{R}^{d}} S\bigg(\frac{\nu^{\mathrm{ac}}_{N}(x)}{\det(\nabla^{2} \phi_j(x))}\bigg) \ud \nu_{N} = \mathcal{F}_{U, N}(\phi_*; \nu^{\mathrm{ac}}_{N}, \nu_{N}).
\end{equation*}
On the other hand, $\phi_{j} \to \phi_\ast$ in $C^{2, \beta}_{\mathrm{loc}}$ implies $(\nabla \phi_j)_\# \nu_N \to (\nabla \phi_\ast)_\# \nu_N$ narrowly. Then by the lower semi-continuity of $\mathcal{F}_{V}$, $\mathcal{F}_{W}$, and $W_{2}$ under narrow convergence in Lemma~\ref{lemma-7-stability}, we have
\begin{equation*}
    \liminf_{j \to \infty} \JKO_{N}(\phi_j; \nu^{\mathrm{ac}}_{N}, \nu_{N}) \ge \JKO_{N}(\phi_\ast; \nu^{\mathrm{ac}}_{N}, \nu_{N})\,.
\end{equation*}
Since $\phi_j$ is a minimizing sequence, we have 
\begin{equation*}
    \liminf_{j \to \infty} \JKO_{N}(\phi_j; \nu^{\mathrm{ac}}_{N}, \nu_{N}) =  \inf_{\phi \in C_{\Lambda}^{2,\alpha}} \JKO_{N}(\phi; \nu^{\mathrm{ac}}_{N}, \nu_{N})
\end{equation*}
which together with the above inequality implies that $\phi_\ast$ is a minimizer.
\end{proof}

\subsection{$\Gamma$-convergence of single step JKO functional}
In this section, we prove the $\Gamma$-convergence of single-step JKO function. First we recall the following lemmas.
\begin{lemma}\label{lemma-2-1-1}
Let $\mu_N$ and $\mu$ be finite positive Radon measures on
$\mathbb{R}^d$ such that $\mu_N\rightharpoonup\mu$ narrowly.
Suppose that $f\in C_b(\mathbb{R}^d)$, $(f_N)_N \subset C_b(\mathbb{R}^d)$ is bounded and converges locally uniformly to $f$. Then 
\[
\lim_{N\to\infty}\int_{\mathbb{R}^d}f_N\ud\mu_N=
    \int_{\mathbb{R}^d}f\ud\mu.
\]
\end{lemma}
\begin{proof}
This lemma is a classic real analysis result, we place the proof in Appendix~\ref{appendix-a}.
\end{proof}
We now state the  $\Gamma$-convergence result.
\begin{theorem} \label{Gamma-NN}
Given $\nu_N^{\mathrm{ac}},\nu \in P_2(\mathbb{R}^d) \cap C(\mathbb{R}^d)$. Assume that
\begin{equation} \label{1003}
    W_2(\nu_N^{\mathrm{ac}},\nu) \rightarrow 0,
    \quad \text{and} \quad \nu_N^{ac}\to\nu\quad\text{locally uniformly on } \; \mathbb R^d.
\end{equation}
and that 
\begin{equation}\label{10-1-bound}
    -C(1+|x|^2) \le S(\nu_N^{\mathrm{ac}}) \le C
\end{equation}
uniformly in $N$. Let $\{\nu_N\}_{N\ge1}\subset P_2(\mathbb R^d)$ be a sequence of atomic
probability measures such that
\begin{equation*} 
    W_2(\nu_N,\nu) \rightarrow 0\,.
\end{equation*}
Then the functionals $\phi\mapsto
\JKO_N(\phi;\nu_N^{\mathrm{ac}},\nu_N)$
are $\Gamma$-convergent, with respect to local uniform convergence, to $\phi\mapsto
\JKO\bigl((\nabla\phi)_{\#}\nu,\nu\bigr)$.
More precisely, the following properties hold.
\begin{enumerate}
    \item[\rm (i)] $\Gamma\text{-} \textnormal{\textbf{liminf}}$.
    Let $\phi\in C_{\Lambda}^{2,\alpha}(\mathbb R^d)$ and
    let $ \phi_N\in C_{\Lambda}^{2, \alpha}(\mathbb R^d)$ with $\phi_N \to \phi$ locally uniformly.  Then
    \begin{equation*}
        \JKO\bigl((\nabla\phi)_{\#}\nu,\nu\bigr)
        \le
        \liminf_{N\to\infty}
        \JKO_N(\phi_N;\nu_N^{\mathrm{ac}},\nu_N).
    \end{equation*}

    \item[\rm (ii)] $\textnormal{\textbf{Recovery sequence}}$.
    For every
    $\phi\in C_{\Lambda}^{2,\alpha}(\mathbb R^d)$,
    the constant sequence $\phi_N=\phi$ satisfies
    \begin{equation*}
        \lim_{N\to\infty}
        \JKO_N(\phi;\nu_N^{\mathrm{ac}},\nu_N)
        =
        \JKO\bigl((\nabla\phi)_{\#}\nu,\nu\bigr).
    \end{equation*}
\end{enumerate}
\end{theorem}
\begin{proof}We first prove the consistency of replacing integration with empirical sum. Let
\begin{equation*}
    e_{N} := \int_{\mathbb{R}^{d}} S\bigg(\frac{\nu^{\mathrm{ac}}_{N}(x)}{\det\big(\nabla^{2} \phi_{N}(x)\big)}\bigg) \ud \big(\nu_{N} - \nu_{N}^{\mathrm{ac}}\big).
\end{equation*}
Since $\Lambda^{d} \le \det(\nabla^{2} \phi_{N}(x)) \le \Lambda^{-d}$, by~\eqref{10-1-bound}, we have
\begin{equation*}
    \frac{1}{1 + |x|^{2}} S\bigg(\frac{\nu_{N}^{\mathrm{ac}}(x)}{\det\big(\nabla^{2} \phi_{N}(x)\big)}\bigg) \text{ is continuous, bounded} \; \text{and converges locally uniformly to}
\end{equation*}
\begin{equation*}
    \frac{1}{1 + |x|^{2}} S\bigg(\frac{\nu(x)}{\det\big(\nabla^{2} \phi(x)\big)}\bigg)
\end{equation*}
locally uniformly, because $\nu_{N}^{\mathrm{ac}} \to \nu$ uniformly and $\phi_{N} \to \phi$ in $C^{2, \beta}_{\mathrm{loc}}(\mathbb{R}^{d})$ for $0<\beta < \alpha$.

Since $\nu_{N}, \nu_{N}^{\mathrm{ac}} \to \nu \text{ in } P_{2}(\mathbb{R}^{d})$, we have
\begin{equation*}
    (1 + |x|^{2})\nu_{N}^{\mathrm{ac}}, (1 + |x|^{2}) \nu_{N} \text{ are tight sequences of positive finite Radon measures,}
\end{equation*}
which converge to $(1 + |x|^{2}) \nu$ narrowly. Using Lemma~\ref{lemma-2-1-1}, we have
\begin{equation}\label{5-1}
    \begin{aligned}
        \lim_{N \to \infty} e_{N} = \lim_{N \to \infty} \int_{\mathbb{R}^{d}} \frac{1}{1 + |x|^{2}} S\bigg[\frac{\nu^{\mathrm{ac}}_{N}(x)}{\det\big(\nabla^{2} \phi_{N}(x)\big)}\bigg] \times (1 + |x|^{2})\ud \big(\nu_{N} - \nu_{N}^{\mathrm{ac}}\big) = 0.
    \end{aligned}
\end{equation}

For the $\liminf$ inequality, we first consider the internal energy. We have
\begin{equation*}
    \begin{aligned}
        \liminf_{N\to\infty}
        \mathcal{F}_{U,N}
        (\phi_N;\nu_N^{\mathrm{ac}},\nu_N)
        &=
        \liminf_{N\to\infty}
        \left[
        \int_{\mathbb{R}^d}
        S\left(
        \frac{\nu_N^{\mathrm{ac}}(x)}
        {\det(\nabla^2\phi_N(x))}
        \right)
        \,\ud\nu_N^{\mathrm{ac}}(x)
        + e_N
        \right]
        \\
        &=
        \liminf_{N\to\infty}
        \mathcal{F}_{U}
        \big((\nabla\phi_N)_{\#}\nu_N^{\mathrm{ac}}\big)
        \geq
        \mathcal{F}_{U}
        \big((\nabla\phi)_{\#}\nu\big),
    \end{aligned}
\end{equation*}
where the last inequality comes from Lemma~\ref{lemma-7-stability}. Furthermore, the corresponding $\liminf$ inequalities for the confinement energy
$\mathcal{F}_{V}$, the interaction energy $\mathcal{F}_{W}$, and the Wasserstein distance follow directly from the same stability result.
Indeed, whenever $(\nabla\phi_N)_{\#}\nu_N
    \rightharpoonup
    (\nabla\phi)_{\#}\nu$ narrowly, 
we have
\begin{equation*}
        \liminf_{N\to\infty}
        \mathcal{F}_{V}
        \big((\nabla\phi_N)_{\#}\nu_N\big)
        \geq
        \mathcal{F}_{V}
        \big((\nabla\phi)_{\#}\nu\big), \quad
        \liminf_{N\to\infty}
        \mathcal{F}_{W}
        \big((\nabla\phi_N)_{\#}\nu_N\big)
        \geq
        \mathcal{F}_{W}
        \big((\nabla\phi)_{\#}\nu\big),
\end{equation*}
and
\begin{equation*}
    \liminf_{N\to\infty}
    W_2\big((\nabla\phi_N)_{\#}\nu_N,\nu_N\big)
    \geq
    W_2\big((\nabla\phi)_{\#}\nu,\nu\big).
\end{equation*}
Combining these estimates yields the desired $\liminf$ inequality of JKO functional.

For the recovery sequence, we again consider the internal energy first. We decompose the integral as follows:
\begin{equation*}
    \begin{aligned}
        \big|\mathcal{F}_{U, N}(\phi; \nu_{N}^{\mathrm{ac}}, \nu_{N}) - \mathcal{F}_{U}((\nabla \phi) \# \nu)\big| \le \bigg|\int_{\mathbb{R}^{d}} S\bigg[\frac{\nu_{N}^{\mathrm{ac}}(x)}{\det(\nabla^{2} \phi(x))}\bigg] \ud \nu_{N}^{\mathrm{ac}} - \int_{\mathbb{R}^{d}} S\bigg[\frac{\nu(x)}{\det(\nabla^{2} \phi(x))}\bigg] \ud \nu\bigg| + |e_{N}|.
    \end{aligned}
\end{equation*}
Set
$D(x):=\det(\nabla^2\phi(x))$ and define
\begin{equation*}
    g_N(x):=
    \frac{1}{1+|x|^2}
    S\left(\frac{\nu_N^{\mathrm{ac}}(x)}{D(x)}\right),
    \qquad
    g(x):=
    \frac{1}{1+|x|^2}
    S\left(\frac{\nu(x)}{D(x)}\right).
\end{equation*}
then $g_N,g\in C_b(\mathbb R^d)$,
\begin{equation*}
    \sup_N\|g_N\|_{L^\infty(\mathbb R^d)}<\infty,
    \qquad
    g_N\to g
    \quad\text{locally uniformly on }\mathbb R^d.
\end{equation*}
Recall the weighted measures
\begin{equation*}
    \mathrm{d}\widehat\eta_N(x)
    :=(1+|x|^2)\ud\nu_N^{\mathrm{ac}}(x),
    \qquad
    \mathrm{d}\eta(x):=(1+|x|^2)\ud\nu(x).
\end{equation*}
These are finite positive Radon measures because
$\nu_N^{\mathrm{ac}},\nu\in P_2(\mathbb R^d)$.
Moreover, the $W_2$-convergence of $\nu_N^{\mathrm{ac}}$ to $\nu$
and Lemma~\ref{lemma-3-p2} imply that
$\widehat\eta_N\rightharpoonup\eta$ narrowly. Therefore, Lemma~\ref{lemma-2-1-1} yields
\begin{align*}
    \int_{\mathbb R^d}
    S\left(\frac{\nu_N^{\mathrm{ac}}(x)}{D(x)}\right)
    \ud\nu_N^{\mathrm{ac}}(x)
    =
    \int_{\mathbb R^d}g_N\ud\widehat\eta_N
    \to
    \int_{\mathbb R^d}g\ud\eta
    =
    \int_{\mathbb R^d}
    S\left(\frac{\nu(x)}{D(x)}\right)\ud\nu(x)
    =
    \mathcal F_U\bigl((\nabla\phi)_\#\nu\bigr).
\end{align*}
Together with $e_N\to0$ from \eqref{5-1}, applied to the
constant sequence $\phi_N=\phi$, the preceding decomposition
gives
\begin{equation*}
    \mathcal F_{U,N}(\phi;\nu_N^{\mathrm{ac}},\nu_N)
    \to
    \mathcal F_U\bigl((\nabla\phi)_\#\nu\bigr).
\end{equation*}

We next consider the Wasserstein term in the recovery-sequence
construction. Since $\nabla\phi$ is globally Lipschitz, the function
$x\mapsto|\nabla\phi(x)-x|^2$ is continuous and has at most quadratic
growth. Moreover, the convexity of $\phi$ implies that $\nabla\phi$
is an optimal transport map for the quadratic cost. Therefore,
the convergence $\nu_N\to\nu$ in $P_2(\mathbb{R}^d)$ yields
\begin{equation*}
        W_2^2\big((\nabla\phi)_{\#}\nu_N,\nu_N\big)
        =
        \int_{\mathbb{R}^d}
        |\nabla\phi(x)-x|^2\,\ud\nu_N(x)
        \to
        \int_{\mathbb{R}^d}
        |\nabla\phi(x)-x|^2\,\ud\nu(x)
        =
        W_2^2\big((\nabla\phi)_{\#}\nu,\nu\big).
\end{equation*}
Finally, since $V$ and $W$ are assumed to be sub-quadratic and $(\nabla\phi)_{\#}\nu_N \to (\nabla\phi)_{\#}\nu$ in $W_{2}$ metric, we have
\begin{equation*}
        \lim_{N\to\infty}
        \mathcal{F}_V\bigl((\nabla\phi)_{\#}\nu_N\bigr)
        =
        \mathcal{F}_V\bigl((\nabla\phi)_{\#}\nu\bigr),\quad
        \lim_{N\to\infty}
        \mathcal{F}_W\bigl((\nabla\phi)_{\#}\nu_N\bigr)
        =
        \mathcal{F}_W\bigl((\nabla\phi)_{\#}\nu\bigr).
\end{equation*}
by the continuity of the confinement and interaction energies under the $W_{2}$ convergence.

This completes the recovery-sequence construction and, together with the $\liminf$ inequality, proves the $\Gamma$-convergence of the JKO functionals.
\end{proof}

\subsection{Convergence in the large sample limit}
In this part, we show the convergence of the particle scheme based on the above $\Gamma$-convergence result.

We establish the following lemma as the result of Proposition~\ref{theorem-10-concentration}, which states that regularity requirements of the $\Gamma$-convergence can be iterated along $\tk = 1, \cdots, \infty$.

\begin{lemma}[iteration]\label{lemma-18-iteration} Let $\mu_{\tk, N}, \mu_{\tk, N}^{\mathrm{ac}} \in P_{2}(\mathbb{R}^{d})$ be the solution of JKO scheme at $\tk$ step such that
\begin{equation*}
    -C_{\tk}(1 + |x|^{2}) \le S\big(\mu_{\tk, N}^{\mathrm{ac}}(x)\big) \le C_{\tk}\,,
\end{equation*}
and let $\phi_{\tk + 1, N}$ be the Brenier potential solution at $\tk+1$-step (i.e., solution to \eqref{1002} with $\nu_N$ and $\nu_N^\text{ac}$ replaced by $\mu_N$ and $\mu_N^\text{ac}$ ), and let
\begin{equation*}
    (\mu_{\tk + 1, N}, \mu_{\tk + 1, N}^{\mathrm{ac}}) = {\nabla \phi_{\tk + 1, N}}_\# (\mu_{\tk, N}, \mu_{\tk, N}^{\mathrm{ac}})\,.
\end{equation*}
Then we have
\begin{equation*}
    |x_{0}(\phi_{\tk + 1, N})| \le R_{\tk},\quad -C_{\tk + 1}(1 + |x|^{2}) \le S(\mu_{\tk + 1, N}^{\mathrm{ac}}(x)) \le C_{\tk + 1}
\end{equation*}
for some constant $C_{\tk + 1}$. The former condition implies that $\mu_{\tk + 1, N}, \mu_{\tk + 1, N}^{\mathrm{ac}} \in P_{2}(\mathbb{R}^{d})$.
\end{lemma}

\begin{proof}[Proof of Theorem~\ref{theorem-3-neural}]
We prove the result by induction. Let $\phi_{\tk}$ and
$\phi_{\tk,N}$ be the continuum and discrete minimizing potentials,
respectively. At each fixed step, the admissible class
$C_{\Lambda}^{2,\alpha}(\mathbb R^d)$ contains $\phi_{\tk}$,
with all its defining constants independent of $N$; their possible
dependence on $\tk$ is suppressed.

Assume that $\mu_{\tk-1,N},\mu_{\tk-1,N}^{\mathrm{ac}}\to
\mu_{\tk-1}$ in $W_2$, that the auxiliary densities converge
locally uniformly, and that they satisfy a bound of the form~\eqref{5-1-bound}
with a constant independent of $N$. The densities are continuous
and strictly positive. These assertions hold initially since
$\mu_{0,N}^{\mathrm{ac}}=\mu_0$ and by the assumptions on the
initial data.

Theorem~\ref{theorem-1002-existence} gives the discrete minimizers.
Since the input second moments are uniformly bounded,
Lemma~\ref{lemma-18-iteration} gives
$\sup_N|x_0(\phi_{\tk,N})|<\infty$.
Normalize $\phi_{\tk,N}(0)=\phi_{\tk}(0)=0$.
By Lemma~\ref{lemma-18-compact}, a subsequence, not relabeled,
converges to $\phi_{*,\tk}\in C_{\Lambda}^{2,\alpha}(\mathbb R^d)$:
\begin{equation*}
    \phi_{\tk,N}\rightarrow\phi_{*,\tk}
    \quad\text{in }C_{\mathrm{loc}}^{2,\beta}(\mathbb R^d),
    \qquad 0<\beta<\alpha.
\end{equation*}
Continuum minimality and the $\Gamma$-liminf inequality in
Theorem~\ref{Gamma-NN}, with locally uniform density convergence
in~\eqref{1003}, give 
\begin{align*}
    \JKO
    ((\nabla\phi_{\tk})_\#\mu_{\tk-1},\mu_{\tk-1})
    \le
    \JKO
    ((\nabla\phi_{*,\tk})_\#\mu_{\tk-1},\mu_{\tk-1})\le
    \liminf_{N\to\infty}
    \JKO_N
    (\phi_{\tk,N};\mu_{\tk-1,N}^{\mathrm{ac}},\mu_{\tk-1,N}).
\end{align*}
On the other hand, discrete minimality and the constant recovery
sequence $\phi_{\tk}$ give 
\begin{align*}
    \limsup_{N\to\infty}
    \JKO_N
    (\phi_{\tk,N};\mu_{\tk-1,N}^{\mathrm{ac}},\mu_{\tk-1,N})\le
    \lim_{N\to\infty}
    \JKO_N
    (\phi_{\tk};\mu_{\tk-1,N}^{\mathrm{ac}},\mu_{\tk-1,N})
    =\JKO
    ((\nabla\phi_{\tk})_\#\mu_{\tk-1},\mu_{\tk-1}).
\end{align*}
Thus $(\nabla\phi_{*,\tk})_\#\mu_{\tk-1}=\mu_{\tk}$ by
uniqueness of the continuum JKO minimizer. Brenier's uniqueness
theorem, strict positivity of $\mu_{\tk-1}$, and continuity of
the gradients imply
$\nabla\phi_{*,\tk}=\nabla\phi_{\tk}$ everywhere.
The normalization then gives $\phi_{*,\tk}=\phi_{\tk}$.
Relative compactness and uniqueness of the subsequential limit
therefore yield convergence of the entire sequence:
\begin{equation*}
    \phi_{\tk,N}\to\phi_{\tk}
    \quad\text{in }C_{\mathrm{loc}}^{2,\beta}(\mathbb R^d),
    \qquad 0<\beta<\alpha.
\end{equation*}

We next prove convergence of the measures in $W_2$.
The uniform Lipschitz bounds on the gradients give
\begin{equation*}
    \begin{aligned}
        W_2(\mu_{k,N},\mu_k)
        &=
        W_2\bigl(
            (\nabla\phi_{k,N})_{\#}\mu_{k-1,N},
            (\nabla\phi_k)_{\#}\mu_{k-1}
        \bigr)\\
        &\leq
        \Lambda^{-1}W_2(\mu_{k-1,N},\mu_{k-1})
        +
        \left(
            \int_{\mathbb{R}^d}
            |\nabla\phi_{k,N}(x)-\nabla\phi_k(x)|^2
            \,\ud\mu_{k-1}(x)
        \right)^{1/2}.
    \end{aligned}
\end{equation*}
The first term tends to zero by induction. The integrand tends
to zero pointwise and is bounded by $A_{\tk}(1+|x|^2)$,
with $A_{\tk}$ independent of $N$, by convergence at the origin
and the uniform Lipschitz bounds. Dominated convergence therefore
shows that the second term tends to zero. Applying the same
estimate with $\mu_{\tk-1,N}^{\mathrm{ac}}$ in place of
$\mu_{\tk-1,N}$ also gives
\begin{equation*}
    \mu_{\tk,N}\to\mu_{\tk},
    \qquad
    \mu_{\tk,N}^{\mathrm{ac}}\to\mu_{\tk}
    \quad\text{in }W_2.
\end{equation*}

Finally, as in Lemma~\ref{lemma-18-iteration}, the gradient maps
are global $C^1$ diffeomorphisms with $\Lambda^{-1}$-Lipschitz
inverses. Set
$Z_N:=(\nabla\phi_{\tk,N})^{-1}$ and
$Z:=(\nabla\phi_{\tk})^{-1}$. Since
\begin{equation*}
    |Z_N(y)-Z(y)|
    \le\Lambda^{-1}
    |\nabla\phi_{\tk,N}(Z(y))-\nabla\phi_{\tk}(Z(y))|,
\end{equation*}
we have $Z_N\to Z$ locally uniformly. The change-of-variables
formula yields
\begin{align*}
    \mu_{\tk,N}^{\mathrm{ac}}(y)
    =\frac{\mu_{\tk-1,N}^{\mathrm{ac}}(Z_N(y))}
    {\det\bigl(\nabla^2\phi_{\tk,N}(Z_N(y))\bigr)}\to
    \frac{\mu_{\tk-1}(Z(y))}
    {\det\bigl(\nabla^2\phi_{\tk}(Z(y))\bigr)}
    =\mu_{\tk}(y)
    \quad\text{locally uniformly}.
\end{align*}
Indeed, the inverse images of each fixed ball lie in a common
compact set for all sufficiently large $N$. On that set, the
previous densities converge locally uniformly by induction,
the Hessians converge uniformly, and their determinants are
bounded below by $\Lambda^d$.
The formula also preserves continuity and strict positivity.
Together with the uniform density bound supplied by
Lemma~\ref{lemma-18-iteration}, this verifies the hypotheses of
Theorem~\ref{Gamma-NN} at the next step and completes the induction.
\end{proof}

\begin{proof}[Proof of Lemma~\ref{lemma-18-iteration}]
Fix $\tk$. By the induction hypothesis, the input measures
have uniformly bounded second moments:
\begin{equation*}
    \sup_N\left(
        M_2(\mu_{\tk,N})
        +M_2(\mu_{\tk,N}^{\mathrm{ac}})
    \right)<\infty.
\end{equation*}
The identity-map potential $q(x)=|x|^2/2$ belongs to the
admissible class. By minimality, the bound
$S(\mu_{\tk,N}^{\mathrm{ac}})\le C_{\tk}$, and the
subquadratic growth of $V$ and $W$, we have
\begin{align*}
    \JKO_N
    (\phi_{\tk+1,N};
     \mu_{\tk,N}^{\mathrm{ac}},\mu_{\tk,N})
    &\le
    \JKO_N
    (q;\mu_{\tk,N}^{\mathrm{ac}},\mu_{\tk,N})\\
    &=
    \int_{\mathbb R^d}
    S(\mu_{\tk,N}^{\mathrm{ac}}(x))\ud\mu_{\tk,N}(x)
    +\mathcal F_V(\mu_{\tk,N})+ \mathcal F_W(\mu_{\tk,N}) \le E_{\tk},
\end{align*}
where $E_{\tk}$ is independent of $N$.
The coercivity estimate in the proof of
Proposition~\ref{theorem-10-concentration}, together with
the uniform input-moment bound, therefore yields
\begin{equation*}
    |x_0(\phi_{\tk+1,N})|\le R_{\tk}
\end{equation*}
for some $R_{\tk}$ independent of $N$.

We next propagate the density bound in~\eqref{5-1-bound}.
Set
\begin{equation*}
    D_N(x):=\det(\nabla^2\phi_{\tk+1,N}(x)).
\end{equation*}
The uniform Hessian bounds imply that
$\nabla\phi_{\tk+1,N}$ is a global $C^1$ diffeomorphism.
Indeed, for every $y\in\mathbb R^d$, the function
$x\mapsto\phi_{\tk+1,N}(x)-y\cdot x$ is strongly convex
and coercive, so it has a unique minimizer satisfying
$\nabla\phi_{\tk+1,N}(x)=y$. The Hessian is everywhere
invertible, so the inverse function theorem applies.

Since
\begin{equation*}
    \mu_{\tk+1,N}^{\mathrm{ac}}
    =
    (\nabla\phi_{\tk+1,N})_\#
    \mu_{\tk,N}^{\mathrm{ac}},
\end{equation*}
the change-of-variables formula gives
\begin{equation*}
    \mu_{\tk+1,N}^{\mathrm{ac}}
    \bigl(\nabla\phi_{\tk+1,N}(x)\bigr)
    =
    \frac{\mu_{\tk,N}^{\mathrm{ac}}(x)}{D_N(x)}.
\end{equation*}
Consequently, the explicit form of $S$ yields
\begin{equation*}
    S\left(
        \mu_{\tk+1,N}^{\mathrm{ac}}
        \bigl(\nabla\phi_{\tk+1,N}(x)\bigr)
    \right)
    =
    \begin{cases}
        S(\mu_{\tk,N}^{\mathrm{ac}}(x))
        -\log D_N(x), & m=1,\\[1mm]
        D_N(x)^{1-m}
        S(\mu_{\tk,N}^{\mathrm{ac}}(x)), & m\ne1.
    \end{cases}
\end{equation*}
For $m>1$, we use the continuous extension $S(0)=0$.

Using
\begin{equation*}
    \Lambda^d\le D_N(x)\le\Lambda^{-d}
    \qquad\text{and}\qquad
    -C_{\tk}(1+|x|^2)
    \le S(\mu_{\tk,N}^{\mathrm{ac}}(x))
    \le C_{\tk},
\end{equation*}
we obtain
\begin{equation*}
    -A_{\tk}(1+|x|^2)
    \le
    S\left(
        \mu_{\tk+1,N}^{\mathrm{ac}}
        \bigl(\nabla\phi_{\tk+1,N}(x)\bigr)
    \right)
    \le A_{\tk},
\end{equation*}
where $A_{\tk}$ is independent of $N$. For example, one
may take
\begin{equation*}
    A_{\tk}:=
    \begin{cases}
        C_{\tk}+d|\log\Lambda|, & m=1,\\
        \Lambda^{-d|1-m|}C_{\tk}, & m\ne1.
    \end{cases}
\end{equation*}

To express this estimate in the destination coordinate,
set $y=\nabla\phi_{\tk+1,N}(x)$.
By Lemma~\ref{lemma:15-equi},
\begin{equation*}
    \Lambda
    |x-x_0(\phi_{\tk+1,N})|
    \le |\nabla\phi_{\tk+1,N}(x)|=|y|.
\end{equation*}
Hence
\begin{equation*}
    |x|\le R_{\tk}+\Lambda^{-1}|y|,
    \qquad
    1+|x|^2
    \le 1+2R_{\tk}^2+2\Lambda^{-2}|y|^2
    \le B_{\tk}(1+|y|^2),
\end{equation*}
where
\begin{equation*}
    B_{\tk}:=\max\{1+2R_{\tk}^2,\,2\Lambda^{-2}\}.
\end{equation*}
Since $\nabla\phi_{\tk+1,N}$ is onto $\mathbb R^d$,
we conclude that
\begin{equation*}
    -C_{\tk+1}(1+|y|^2)
    \le S(\mu_{\tk+1,N}^{\mathrm{ac}}(y))
    \le C_{\tk+1}
    \qquad\text{for every }y\in\mathbb R^d,
\end{equation*}
where $C_{\tk+1}:=A_{\tk}B_{\tk}$ is independent of $N$.

Finally, the upper bound in Lemma~\ref{lemma:15-equi}
gives
\begin{equation*}
    |\nabla\phi_{\tk+1,N}(x)|^2
    \le
    \Lambda^{-2}|x-x_0(\phi_{\tk+1,N})|^2
    \le
    2\Lambda^{-2}(|x|^2+R_{\tk}^2).
\end{equation*}
Integrating against the two input measures yields 
\begin{align*}
    M_2(\mu_{\tk+1,N})\le
    2\Lambda^{-2}
    \bigl(M_2(\mu_{\tk,N})+R_{\tk}^2\bigr), \qquad 
    M_2(\mu_{\tk+1,N}^{\mathrm{ac}})\le
    2\Lambda^{-2}
    \bigl(M_2(\mu_{\tk,N}^{\mathrm{ac}})+R_{\tk}^2\bigr).
\end{align*}
Thus both output measures belong to $ P_2(\mathbb R^d)$,
with second moments uniformly bounded in $N$.
This completes the proof.
\end{proof}

\subsection*{Acknowledgments} 
Part of this work is supported by the National Science Foundation under Grant No. DMS-
2424139, when LW and WG visited the Simons Laufer Mathematical Sciences Institute in Berkeley,
California, during the Fall 2025 semester to participate the semester program on “Kinetic Theory: Novel
Statistical, Stochastic and Analytical Methods".
The research of WG was additionally supported by NSF grant DMS-2154578. LW was additionally supported by NSF grant DMS-2513336 and Simons Foundation TSM. The authors would like to thank Professors Jos\'e Carrillo and Matt Jacobs for fruitful discussions. 

~\\
\textbf{AI disclosure}: GPT Pro was used to assist with drafting, editing, and verification and revision of minor technical details. It was not used in developing the main mathematical ideas or core contributions of the work.

\appendix
\section{Proof of technical lemmas}\label{appendix-a}
\begin{proof}[Proof of Lemma~\ref{lemma-2-1-1}]
Narrow convergence, tested against the constant function $1$,
implies that
\begin{equation*}
    K:=\sup_N\mu_N(\mathbb{R}^d)<\infty.
\end{equation*}
Let $M$ be a common bound for $\|f\|_{L^\infty}$ and
$\|f_N\|_{L^\infty}$.

For $R>0$, choose $\chi_R\in C_c(\mathbb{R}^d)$ such that
$0\le\chi_R\le1$, $\chi_R=1$ on $B_R(0)$, and
$\operatorname{supp}\chi_R\subset B_{2R}(0)$. Then 
\begin{align*}
\left|
\int_{\mathbb{R}^d}f_N\ud\mu_N
-\int_{\mathbb{R}^d}f\ud\mu
\right|
& \le
\int_{\mathbb{R}^d}\chi_R|f_N-f|\ud\mu_N
+\int_{\mathbb{R}^d}(1-\chi_R)|f_N-f|\ud\mu_N+
\left|
\int_{\mathbb{R}^d}f\ud\mu_N
-\int_{\mathbb{R}^d}f\ud\mu
\right| \\
& \le
K\sup_{B_{2R}(0)}|f_N-f|
+2M\int_{\mathbb{R}^d}(1-\chi_R)\ud\mu_N+
\left|
\int_{\mathbb{R}^d}f\ud\mu_N
-\int_{\mathbb{R}^d}f\ud\mu
\right|.
\end{align*}
For fixed $R$, local uniform convergence and narrow convergence
therefore give
\begin{align*}
\limsup_{N\to\infty}
\left|
\int_{\mathbb{R}^d}f_N\ud\mu_N
-\int_{\mathbb{R}^d}f\ud\mu
\right|\le
2M\int_{\mathbb{R}^d}(1-\chi_R)\ud\mu \le
2M\mu\bigl(\mathbb{R}^d\setminus B_R(0)\bigr).
\end{align*}
Letting $R\to\infty$ proves the claim, since $\mu$ is finite.
\end{proof}

\bibliographystyle{alpha}
\bibliography{JKO_conv}

@book{ambrosio2008gradient,
  title={Gradient flows: in metric spaces and in the space of probability measures},
  author={Ambrosio, Luigi and Gigli, Nicola and Savar{\'e}, Giuseppe},
  year={2008},
  publisher={Springer Science \& Business Media}
}

@article{JKO,
    author = {Jordan, Richard and Kinderlehrer, David and Otto, Felix},
    title = {The Variational Formulation of the Fokker--Planck Equation},
    journal = {SIAM Journal on Mathematical Analysis},
    volume = {29},
    number = {1},
    pages = {1-17},
    year = {1998},
    doi = {10.1137/S0036141096303359},
}

@article{coudreuse2026aronson,
  title={An Aronson-B$\backslash$'enilan/Li-Yau estimate in the JKO scheme in small dimension},
  author={Coudreuse, Fanch},
  journal={arXiv preprint arXiv:2604.04169},
  year={2026}
}

@book{santambrogio2015optimal,
  title={Optimal transport for applied mathematicians},
  author={Santambrogio, Filippo},
  volume={87},
  year={2015},
  publisher={Springer}
}

@article{JungeMatthesOsberger2017,
  author    = {Junge, Oliver and Matthes, Daniel and Osberger, Horst},
  title     = {A fully discrete variational scheme for solving nonlinear {F}okker--{P}lanck equations in multiple space dimensions},
  journal   = {SIAM Journal on Numerical Analysis},
  volume    = {55},
  number    = {1},
  pages     = {419--443},
  year      = {2017},
  publisher = {Society for Industrial and Applied Mathematics},
  doi       = {10.1137/16M1056560},
  issn      = {0036-1429}
}

@article{GangboMcCann1996,
  author  = {Gangbo, Wilfrid and McCann, Robert J.},
  title   = {The geometry of optimal transportation},
  journal = {Acta Mathematica},
  year    = {1996},
  volume  = {177},
  number  = {2},
  pages   = {113--161},
  doi     = {10.1007/BF02392620},
  url     = {https://doi.org/10.1007/BF02392620}
}

@article{GangboTudorascu2014,
  author   = {Gangbo, Wilfrid and Tudorascu, Adrian},
  title    = {Weak {KAM} theory on the {Wasserstein}
              torus with multidimensional underlying space},
  journal  = {Comm. Pure Appl. Math.},
  volume   = {67},
  number   = {3},
  year     = {2014},
  pages    = {408--463},
  doi      = {10.1002/cpa.21492},
  mrnumber = {3158572}
  }

@book{villani2008optimal,
  title={Optimal transport: old and new},
  author={Villani, C{\'e}dric},
  volume={338},
  year={2008},
  publisher={Springer}
}

@article{benamou2016discretization,
  title={Discretization of functionals involving the {Monge--Amp{\`e}re} operator},
  author={Benamou, Jean-David and Carlier, Guillaume and M{\'e}rigot, Quentin and Oudet, Edouard},
  journal={Numerische mathematik},
  volume={134},
  number={3},
  pages={611--636},
  year={2016},
  publisher={Springer}
}

@article{otto2001geometry,
  title={The geometry of dissipative evolution equations: the porous medium equation},
  author={Otto, Felix},
  year={2001},
  publisher={Taylor \& Francis},
  journal={},
}

@incollection{carrillo2010particle,
  title={Particle, kinetic, and hydrodynamic models of swarming},
  author={Carrillo, Jos{\'e} A and Fornasier, Massimo and Toscani, Giuseppe and Vecil, Francesco},
  booktitle={Mathematical modeling of collective behavior in socio-economic and life sciences},
  pages={297--336},
  year={2010},
  publisher={Springer}
}

@article{carrillo2019population,
  title={A population dynamics model of cell-cell adhesion incorporating population pressure and density saturation},
  author={Carrillo, Jose A and Murakawa, Hideki and Sato, Makoto and Togashi, Hideru and Trush, Olena},
  journal={Journal of theoretical biology},
  volume={474},
  pages={14--24},
  year={2019},
  publisher={Elsevier}
}

@article{benedetto1997kinetic,
  title={A kinetic equation for granular media},
  author={Benedetto, Dario and Caglioti, Emanuele and Pulvirenti, Mario},
  journal={ESAIM: Mathematical Modelling and Numerical Analysis},
  volume={31},
  number={5},
  pages={615--641},
  year={1997},
  publisher={EDP Sciences}
}

@article{toscani2000one,
  title={One-dimensional kinetic models of granular flows},
  author={Toscani, Giuseppe},
  journal={ESAIM: Mathematical Modelling and Numerical Analysis},
  volume={34},
  number={6},
  pages={1277--1291},
  year={2000},
  publisher={EDP Sciences}
}

@article{keller1970initiation,
  title={Initiation of slime mold aggregation viewed as an instability},
  author={Keller, Evelyn F and Segel, Lee A},
  journal={Journal of theoretical biology},
  volume={26},
  number={3},
  pages={399--415},
  year={1970},
  publisher={Elsevier}
}

@article{daneri2022deterministic,
  title={Deterministic particle approximation of aggregation-diffusion equations on unbounded domains},
  author={Daneri, Sara and Radici, Emanuela and Runa, Eris},
  journal={Journal of Differential Equations},
  volume={312},
  pages={474--517},
  year={2022},
  publisher={Elsevier}
}

@article{di2024optimal,
  title={Optimal transport with nonlinear mobilities: A deterministic particle approximation result},
  author={Di Marino, Simone and Portinale, Lorenzo and Radici, Emanuela},
  journal={Advances in Calculus of Variations},
  volume={17},
  number={3},
  pages={941--974},
  year={2024},
  publisher={De Gruyter}
}

@article{di2015rigorous,
  title={Rigorous derivation of nonlinear scalar conservation laws from follow-the-leader type models via many particle limit},
  author={Di Francesco, Marco and Rosini, Massimiliano D},
  journal={Archive for rational mechanics and analysis},
  volume={217},
  number={3},
  pages={831--871},
  year={2015},
  publisher={Springer}
}

@article{matthes2014convergence,
  title={Convergence of a variational Lagrangian scheme for a nonlineardrift diffusion equation},
  author={Matthes, Daniel and Osberger, Horst},
  journal={ESAIM: Mathematical modelling and numerical analysis},
  volume={48},
  number={3},
  pages={697--726},
  year={2014},
  publisher={EDP Sciences}
}

@article{blanchet2008convergence,
  title={Convergence of the mass-transport steepest descent scheme for the subcritical {Patlak--Keller--Segel} model},
  author={Blanchet, Adrien and Calvez, Vincent and Carrillo, Jos{\'e} A},
  journal={SIAM Journal on Numerical Analysis},
  volume={46},
  number={2},
  pages={691--721},
  year={2008},
  publisher={SIAM}
}

@article{lee2024deep,
  title={Deep JKO: Time-implicit particle methods for general nonlinear gradient flows},
  author={Lee, Wonjun and Wang, Li and Li, Wuchen},
  journal={Journal of Computational Physics},
  volume={514},
  pages={113187},
  year={2024},
  publisher={Elsevier}
}

@article{carrillo2019blob,
  title={A blob method for diffusion},
  author={Carrillo, Jos{\'e} Antonio and Craig, Katy and Patacchini, Francesco S},
  journal={Calculus of Variations and Partial Differential Equations},
  volume={58},
  number={2},
  pages={53},
  year={2019},
  publisher={Springer}
}

@article{carrillo2011global,
  title={Global-in-time weak measure solutions and finite-time aggregation for nonlocal interaction equations},
  author={Carrillo, Jos{\'e} A and DiFrancesco, Marco and Figalli, Alessio and Laurent, Thomas and Slep{\v{c}}ev, Dejan},
  journal={},
  year={2011}
}

@incollection{carrillo2019aggregation,
  title={Aggregation-diffusion equations: dynamics, asymptotics, and singular limits},
  author={Carrillo, Jos{\'e} A and Craig, Katy and Yao, Yao},
  booktitle={Active Particles, Volume 2: Advances in Theory, Models, and Applications},
  pages={65--108},
  year={2019},
  publisher={Springer}
}

@article{kim2012patlak,
  title={The Patlak--Keller--Segel model and its variations: properties of solutions via maximum principle},
  author={Kim, Inwon and Yao, Yao},
  journal={SIAM Journal on Mathematical Analysis},
  volume={44},
  number={2},
  pages={568--602},
  year={2012},
  publisher={SIAM}
}

@article{carrillo2004wasserstein,
  title={{Wasserstein} metric and large-time asymptotics of nonlinear diffusion equations},
  author={Carrillo, Jose Antonio and Toscani, Giuseppe and others},
  journal={New trends in mathematical physics},
  pages={234--244},
  year={2004},
  publisher={World Scientific}
}

@article{carrillo2003kinetic,
  title={Kinetic equilibration rates for granular media and related equations: entropy dissipation and mass transportation estimates},
  author={Carrillo, Jos{\'e} A and McCann, Robert J and Villani, C{\'e}dric},
  journal={Revista Matematica Iberoamericana},
  volume={19},
  number={3},
  pages={971--1018},
  year={2003}
}

@article{carrillo2021lagrangian,
  title={Lagrangian schemes for {Wasserstein} gradient flows},
  author={Carrillo, Jose A and Matthes, Daniel and Wolfram, Marie-Therese},
  journal={Handbook of Numerical Analysis},
  volume={22},
  pages={271--311},
  year={2021},
  publisher={Elsevier}
}

@article{chu2025guaranteeing,
  title={Guaranteeing Higher Order Convergence Rates for Accelerated {Wasserstein} Gradient Flow Schemes},
  author={Chu, Raymond and Jacobs, Matt},
  journal={arXiv preprint arXiv:2511.10884},
  year={2025}
}

@article{carrillo2006contractions,
  title={Contractions in the 2-{Wasserstein} length space and thermalization of granular media},
  author={Carrillo, Jos{\'e} A and McCann, Robert J and Villani, C{\'e}dric},
  journal={Archive for Rational Mechanics and Analysis},
  volume={179},
  number={2},
  pages={217--263},
  year={2006},
  publisher={Springer}
}

@article{westdickenberg2010variational,
  title={Variational particle schemes for the porous medium equation and for the system of isentropic Euler equations},
  author={Westdickenberg, Michael and Wilkening, Jon},
  journal={ESAIM: Mathematical Modelling and Numerical Analysis},
  volume={44},
  number={1},
  pages={133--166},
  year={2010},
  publisher={EDP Sciences}
}

@article{craig2023blob,
  title={A blob method for inhomogeneous diffusion with applications to multi-agent control and sampling},
  author={Craig, Katy and Elamvazhuthi, Karthik and Haberland, Matt and Turanova, Olga},
  journal={Mathematics of Computation},
  volume={92},
  number={344},
  pages={2575--2654},
  year={2023}
}

@article{matthes2019variational,
  title={A variational formulation of the BDF2 method for metric gradient flows},
  author={Matthes, Daniel and Plazotta, Simon},
  journal={ESAIM: Mathematical Modelling and Numerical Analysis},
  volume={53},
  number={1},
  pages={145--172},
  year={2019},
  publisher={EDP Sciences}
}

@article{feng2026learn,
  title={Learn to Evolve: Self-supervised Neural JKO Operator for {Wasserstein} Gradient Flow},
  author={Feng, Xue and Wang, Li and Needell, Deanna and Lai, Rongjie},
  journal={arXiv preprint arXiv:2601.05583},
  year={2026}
}

@article{craig2016blob,
  title={A blob method for the aggregation equation},
  author={Craig, Katy and Bertozzi, Andrea},
  journal={Mathematics of computation},
  volume={85},
  number={300},
  pages={1681--1717},
  year={2016}
}

@article{carrillo2022primal,
  title={Primal dual methods for Wasserstein gradient flows},
  author={Carrillo, Jos{\'e} A and Craig, Katy and Wang, Li and Wei, Chaozhen},
  journal={Foundations of Computational Mathematics},
  volume={22},
  number={2},
  pages={389--443},
  year={2022},
  publisher={Springer}
}

@article{liu2017random,
  title={A random particle blob method for the Keller-Segel equation and convergence analysis},
  author={Liu, Jian-Guo and Yang, Rong},
  journal={Mathematics of Computation},
  volume={86},
  number={304},
  pages={725--745},
  year={2017}
}

@article{craig2016convergence,
  title={Convergence of regularized nonlocal interaction energies},
  author={Craig, Katy and Topaloglu, Ihsan},
  journal={SIAM Journal on Mathematical Analysis},
  volume={48},
  number={1},
  pages={34--60},
  year={2016},
  publisher={SIAM}
}

@article{carrillo2016convergence,
  title={Convergence of a particle method for diffusive gradient flows in one dimension},
  author={Carrillo, Jos{\'e} A and Patacchini, Francesco Saverio and Sternberg, Peter and Wolansky, Gershon},
  journal={SIAM Journal on Mathematical Analysis},
  volume={48},
  number={6},
  pages={3708--3741},
  year={2016},
  publisher={SIAM}
}

@book{gilbarg1977elliptic,
  title={Elliptic partial differential equations of second order},
  author={Gilbarg, David and Trudinger, Neil S and Gilbarg, David and Trudinger, NS},
  volume={224},
  number={2},
  year={1977},
  publisher={Springer}
}

@book{braides2002gamma,
    title={Gamma-convergence for Beginners},
    author={Braides, Andrea},
    volume={22},
    year={2002},
    publisher={Clarendon Press}
}

@book{braides2014local,
    title={Local minimization, variational evolution and $\Gamma$-convergence},
    author={Braides, Andrea},
    volume={2094},
    year={2014},
    publisher={Springer}
}

@article{korotin2022neural,
  title={Neural optimal transport},
  author={Korotin, Alexander and Selikhanovych, Daniil and Burnaev, Evgeny},
  journal={arXiv preprint arXiv:2201.12220},
  year={2022}
}

@inproceedings{makkuva2020optimal,
  title={Optimal transport mapping via input convex neural networks},
  author={Makkuva, Ashok and Taghvaei, Amirhossein and Oh, Sewoong and Lee, Jason},
  booktitle={International Conference on Machine Learning},
  pages={6672--6681},
  year={2020},
  organization={PMLR}
}

@article{han2023high,
  title={High order schemes for gradient flow with respect to a metric},
  author={Han, Saem and Esedoḡlu, Selim and Garikipati, Krishna},
  journal={Journal of Computational Physics},
  volume={494},
  pages={112516},
  year={2023},
  publisher={Elsevier}
}

@article{dumont2026learning,
  title={Learning {Monge} maps by lifting and constraining {Wasserstein} gradient flows},
  author={Dumont, Th{\'e}o and Lacombe, Th{\'e}o and Vialard, Fran{\c{c}}ois-Xavier},
  journal={},
  year={2026}
}

@article{coudreuse2025li,
  title={Li-Yau-Hamilton Inequality on the JKO Scheme for the Granular-Medium Equation},
  author={Coudreuse, Fanch},
  journal={arXiv preprint arXiv:2510.09231},
  year={2025}
}

@article{brenier1987decomposition,
  title={D{\'e}composition polaire et r{\'e}arrangement monotone des champs de vecteurs},
  author={Brenier, Yann},
  journal={CR Acad. Sci. Paris S{\'e}r. I Math.},
  volume={305},
  pages={805--808},
  year={1987}
}

\end{document}